\documentclass[12pt,letterpaper,reqno]{amsart}

\usepackage{graphicx}
\usepackage{epstopdf}
\usepackage{fullpage}
\usepackage{amsmath}
\usepackage{amsfonts}
\usepackage{amssymb}
\usepackage{mathrsfs}
\usepackage{dsfont}
\usepackage{xargs}

\newtheorem{thm}{Theorem}[section]
\newtheorem{prop}[thm]{Proposition}
\newtheorem{lemma}[thm]{Lemma}

\theoremstyle{definition}

\newtheorem{defn}[thm]{Definition}
\newtheorem{rmk}[thm]{Remark}
\newtheorem*{prf}{Proof}
\numberwithin{equation}{section}
\numberwithin{figure}{section}
\usepackage{theoremref}

\newcommand{\floor}[1]{\lfloor #1 \rfloor}

\newcommand{\hs}[1][0.5cm]{\hspace{#1}}

\newcommand{\R}[0]{\mathbb{R}}
\newcommand{\C}[0]{\mathbb{C}}
\newcommand{\Z}[0]{\mathbb{Z}}
\newcommand{\N}[0]{\mathbb{N}}

\newcommand{\1}[0]{\mathds{1}}
\newcommand{\vphi}[0]{\varphi}

\begin{document}
\date{\today}
\title{Constructing Tight Gabor Frames Using CAZAC sequences}
\author{Mark Magsino}

\begin{abstract}
The construction of finite tight Gabor frames plays an important role in
many applications.
These applications include significant ones in signal and image processing. 
We explore when constant amplitude zero autocorrelation (CAZAC) sequences can be used to generate
tight Gabor frames.
The main theorem uses
Janssen's representation and the zeros of the discrete periodic ambiguity function to give 
necessary and sufficient conditions for determining whether any Gabor frame is tight.
The relevance of the theorem depends significantly on the construction of examples.
These examples are necessarily intricate, and to a large extent, depend on CAZAC sequences.
Finally, we present an alternative method for determining whether
a Gabor system yields a tight frame. This alternative method does not prove tightness
using the main theorem, but instead uses the Gram matrix of the Gabor system.
\end{abstract}

\maketitle 

\section{Introduction}
\subsection{Background}


Frames were introduced in 1952 by Duffin and Schaeffer \cite{DufSch1952} in their research on
nonharmonic Fourier series. They used frames to compute the coefficients of a linear combination of 
vectors which were linearly dependent and spanned its Hilbert space. Since then, frames have 
been used in applications such as the analysis of wavelets, and in signal and image
processing \cite{Dau1992}\cite{Mal1999}\cite{Mix2012}\cite{PfaWal2006}\cite{Str2001}. 
Frames can be viewed as a generalization
of orthonormal bases for Hilbert spaces. Like bases, frames still span the Hilbert space, 
but unlike bases they are allowed to be
linearly dependent. In the context of signal processing, the primary advantage of frames is that 
they provide stable representations of signals which are robust in the presense of erasures and 
noise \cite{KovCev2007}\cite{KovCev2007b}. Frames for finite vector spaces, i.e., finite frames, 
are of particular interest for engineering or computational applications, and as such there 
has been significant research conducted on finite frames 
\cite{BalCasHeiLan2006}\cite{BalCasHeiLan2006b}\cite{CasKut2013}\cite{Chr2008}.


Let $\mathcal{H}$ be a separable Hilbert space over the complex field $\C$. A sequence
$\mathcal{F} = \{v_i\} \subseteq \mathcal{H}$ is a \textit{frame} for $\mathcal{H}$ if there 
exist $A,B > 0$ such that
\begin{equation}\label{framedefn}
    \forall x \in \mathcal{H}, \text{ } 
    A\|x\|_2^2 \leq \sum_{v_i \in \mathcal{F}} |\langle x, v_i \rangle|^2 \leq B\|x\|_2^2.
\end{equation}
Since we want to view frames as a 
generalization of orthonormal bases, we want to be able to write any vector $x \in \mathcal{H}$ 
in terms of $v_i \in \mathcal{F}$.
If $\mathcal{F}$ is a frame, then we can write any $x \in \mathcal{H}$ as the linear combination,
\begin{equation}\label{framereconst}
    x = \sum_{v_i \in \mathcal{F}} \langle x, v_i \rangle S^{-1} v_i,
\end{equation}
where $S$ is a well-defined linear operator associated with $\mathcal{F}$ known as the 
\textit{frame operator} of $\mathcal{F}$.
In general, it is non-trivial to compute the invese of the frame operator. However, if is
possible to have $A = B$, then we have the special case of a \textit{tight frame}.
In this case, $S = A \, I\!d$, and (\ref{framereconst}) can be re-written as
\begin{equation}\label{tightreconst}
    x = \frac{1}{A} \sum_{v_i \in \mathcal{F}} \langle x, v_i \rangle v_i.
\end{equation}
This makes tight frames particularly desirable since (\ref{tightreconst}) is computationally
easier than (\ref{framereconst}). This has motivated research into the discovery and construction 
of tight frames 
\cite{BenFic2003}\cite{CasKov2003}\cite{DauHanRonShe2003}\cite{ValWal2004}\cite{ValWal2008},
as well as the transformation of frames into tight frames \cite{KutOkoPhiTul2013}.


Finite frames are sometimes studied in the context of time-frequency analysis.
The beginnings of time-frequency analysis go back to Gabor's
paper on communication theory \cite{Gab1946}, where he used time-frequency rectangles  
to simultaenously analyze the time and frequency content of Gaussian functions.
These methods are restricted by the Heisenberg uncertainty principle, which essentially states
that no function can be simultaneously well concentrated in both time and frequency 
\cite{Ben1993}\cite{FolSit1997}.
It is this limitation that makes the mathematical theory of time-frequency analysis
both difficult and interesting, and there is significant research in the area of time-frequency
analysis
\cite{BenHeiWal1998}\cite{Dau1990}\cite{LawPfaWal2005}\cite{PfaRau2010}\cite{PfaRauTro2013}. 

In particular, systems which consist of translations and modulations of a generating
function are called Gabor systems. To ensure any signal (function) can be constructed from a Gabor
system, one would like to show that the system forms an orthonormal basis, or more generally,
a frame. Gabor's original suggestion was to use Gaussian functions as the generator. 
His suggestion did not generate a Riesz basis for $L^2(\R)$ \cite{Jan1981}, 
the space of square integrable functions, but a minor alteration to his suggestion does make the
system into a frame for $L^2(\R)$ \cite{Lyu1992}. This example motivates studying when Gabor
systems generated by other functions are tight frames. 
An excellent reference and exposition on time-frequency
analysis is \cite{Gro2013}. 


Constant amplitude zero autocorrelation (CAZAC) sequences are often used in radar and 
communication theory for various applications \cite{LevMoz2004} such as
error-correcting codes. The zero autocorrelation property gives that CAZAC sequences and their
translates are orthogonal. One might hope similar sparsity properties extend off the zero time
axis. This is the motivation which suggests that CAZAC sequences may be suitable for generating
tight Gabor frames.


\subsection{Theme}

The central idea is to use CAZAC sequences to generate tight Gabor frames. The motivation
is based on the following train of thought. Frames are useful for computational purposes 
since they are
often stable representations which are robust in the presence of erasures and noise. Tight
frames are even more computationally convenient because reconstruction avoids the need of
computing the inverse of the frame operator. Finite tight frames are then considered since
computers can ultimately only do finite computations. Finite tight Gabor frames are used to
see the role of time and frequency in signals we would like to represent with tight frames.
Finally, CAZACs are used in radar and communication theory, and so it is natural to consider
the role of CAZACs in the context of Gabor frames. 

One of the main themes is constructing Gabor frames using subgroups of the time-frequency group.
The idea is that tight Gabor frames using subgroups would have fewer coefficients and require less
computation while maintaining group structure. In this context, our optimal theoretical
result is \thref{aptightframe}.
Because of the intricacies involved in quantifying \thref{aptightframe},
we construct nontrivial examples of \thref{aptightframe}. 

Of comparable relevance are CAZAC sequences. 
\thref{aptightframe} requires suitable sparsity in the ambiguity function of the generating
sequence.
Several CAZAC sequences can be shown to have the required sparsity in their amibguity function
in several cases.
As such, we leverage this sparsity and use these CAZAC sequences to construct nontrivial
examples of \thref{aptightframe}.
This prompts the question of whether
CAZACs are always suitable for constructing tight Gabor frames, and if not, which ones are
suitable.


\subsection{Outline}

In Section \ref{sectfdpaf} we give a necessary and sufficient condition for proving a Gabor
system in $\C^N$ is a tight frame. This is accomplished through Janssen's representation and
by showing sufficient sparsity in the discrete periodic ambiguity function. 
In Section \ref{cazacsec}, we give a brief overview on constant amplitude zero autocorrelation
(CAZAC) sequences and define the CAZAC sequences used in the examples:
Chu, P4, Wiener, square length Bj\"orck-Saffari, and Milewski sequences. 
We also provide two alternative
formulations of the question of constructing new CAZAC sequences.

In Section \ref{cazacdpaf}, we compute the discrete periodic ambiguity functions for
the CAZAC sequences we use in the examples. These sequences are the Chu, P4,
Wiener, square length Bj\"orck-Saffari, and Milewski sequences. 
Some details on the Chu, P4, and Wiener sequences can be found in
\cite{BenKonRan2009}. Details on the square length Bj\"orck-Saffari sequence, as well as 
some generalizations of the sequence, can be found in
\cite{BjoSaf1995}. Details on the Milewski sequence and a detailed computation of its
discrete periodic ambiguity function can be found in \cite{BenDon2007} or 
\cite{Mil1983}. 
Section \ref{exsparse} constructs several examples which utilize the sequences from 
Section \ref{cazacdpaf} and \thref{aptightframe}.

We begin Section \ref{gramdpaf} by computing the Gram matrix of a Gabor system in terms
of the discrete periodic ambiguity function of the generating sequence. 
We then use the results of Section \ref{cazacdpaf} in order to write the Gram matrices of 
Gabor systems generated by the Chu, P4, and Wiener sequences.
Section \ref{secgramchu} is focused on an alternative 
method for proving Gabor systems are tight frames. This method proves tightness by
showing the Gram matrix has sufficient rank and that the nonzero eigenvalues of the Gram matrix are the same. Section \ref{secgramchu} begins with the P4 and Chu cases and
is followed up with how to make the necessary adjustments for the Wiener sequence case.

\subsection{Notation}

Let $\vphi \in \C^N$. We denote \textit{translation} by $k \in (\Z/N\Z)$ as $\tau_k$ and 
\textit{modulation} by $\ell \in (\Z/N\Z)$ as
$e_\ell$, and define them by
\begin{equation}\label{thmodform}
    \tau_k(\vphi)[j] = \vphi[j-k] \text{ and } e_\ell(\vphi)[j] = e^{2 \pi i j \ell / N} \vphi[j].
\end{equation}
The \textit{Discrete Fourier Transform} (DFT), $\widehat{\vphi} \in \C^N$, of $\vphi \in \C^N$
is defined by
\begin{equation}\label{dftform}
    \forall \ell \in (\Z/N\Z), \text{ } 
    \widehat{\vphi}[\ell] = \sum_{k = 0}^{N-1} \vphi[k] e^{-2 \pi i k \ell / N}.
\end{equation}
In particular, (\ref{dftform}) is the non-normalized version of the DFT, and so the inverse
is given by
\begin{equation}\label{idftform}
    \vphi[k] = \frac{1}{N} \sum_{\ell = 0}^{N-1} \widehat{\vphi}[\ell] e^{2 \pi i k \ell / N}.
\end{equation}

$\Lambda$ denotes a subgroup of the time-frequency lattice 
$(\Z/N\Z) \times (\Z/N\Z)^{\widehat{}}$, where $(\Z/N\Z)^{\widehat{}} = \{e^{2 \pi i
\ell (\cdot) / N} : \ell \in (\Z/N\Z)\}$ is the group of unimodular characters on $(\Z/N\Z)$.
We choose this notation to emphasize that $(\Z/N\Z)^{\widehat{}}$ is indeed the character group even
though $(\Z/N\Z)^{\widehat{}}$ can be identified with the group $(\Z/N\Z)$. In practice, we shall
use the identification of $(\Z/N\Z)^{\widehat{}}$ with $(\Z/N\Z)$ and
simply write $(k,\ell) \in (\Z/N\Z) \times (\Z/N\Z)^{\widehat{}}$.

\section{Tight frames from sparse discrete periodic ambiguity functions} \label{sectfdpaf}

In this section we give a necessary and sufficient condition for determining when the Gabor
system generated by $\vphi \in \C^N \setminus \{0\}$ and $\Lambda \subseteq (\Z/N\Z) \times 
(\Z/N\Z)^{\widehat{}}$ 
will be a tight frame. The major theme will be: Gabor systems are tight frames if the 
discrete periodic ambiguity functions of their generating functions are sufficently sparse. 
The discrete periodic ambiguity
function is a tool that is often used in radar and communications theory 
\cite{LevMoz2004}, and we will show that
it is closely linked to the short-time Fourier transform.
We will use Janssen's representation to utilize the discrete periodic ambiguity function.
Janssen's representation allows us to write the frame operator as a
linear combination of time-frequency operators whose coefficients can be computed without knowing
the input of the frame operator, cf. Walnut's representation \cite{Pfa2013}\cite{Wal1998}.
To begin, we review finite Gabor systems, the short-time Fourier transform, and the fact that full Gabor systems always generate tight frames. For more details on frame
theory, \cite{Chr2008} is a valuable resource.


Gabor systems in $\C^N$ are families of vectors 
which are generated by a vector $\vphi \in \C^N \setminus \{0\}$ and 
translations and modulations of $\vphi$. Specifically, let $\vphi \in \C^N \setminus \{0\}$ and 
let $\Lambda \subseteq \Z/N\Z \times (\Z/N\Z)^{\widehat{}}$. The family of vectors
\[
(\vphi, \Lambda) = \{ e_\ell \tau_k \vphi : (k,\ell) \in \Lambda\}
\]
is the Gabor system generated by $\vphi$ and $\Lambda$. If there exist $A,B > 0$ such that
the sequence $\mathcal{F} = \{v_i\}_{i=1}^{M} \subseteq \C^N$ satisfy
\begin{equation}\label{framecond}
    \forall x \in \C^N, \text{ } 
    A\|x\|_2^2 \leq \sum_{i=1}^M |\langle x,v_i \rangle|^2 \leq B\|x\|_2^2,
\end{equation}
then $\mathcal{F}$ is said to be a frame for $\C^N$.
Note that $A > 0$ guarantees that $\{v_i\}_{i=1}^M$ spans $\C^N$ and $M \geq N$. 
In fact, $\{v_i\}_{i=1}^M$ spanning $\C^N$ is a necessary and sufficient condition in finite
dimensional Hilbert spaces \cite{Chr2008}.
This allows
us to think of frames as a generalization of bases. If $A = B$ is possible in (\ref{framecond}), 
then $\mathcal{F}$ is 
said to be a tight frame. If a Gabor system $(\vphi, \Lambda)$ satisfies (\ref{framecond}), then 
the Gabor system is said to be a Gabor frame. If, in addition,
$A = B$ is possible in (\ref{framecond}), we call $(\vphi,\Lambda)$ a tight Gabor frame.

Let $\mathcal{F} = \{v_i\}_{i=1}^M$ be a frame for $\C^N$.
The frame operator of $\mathcal{F}$, $S: \C^N \to \C^N$, is defined by
\begin{equation}\label{frameop}
    \forall x \in \C^N, \text{ }S(x) = \sum_{i=1}^M \langle x, v_i \rangle v_i.
\end{equation}
We can reconstruct any vector $x \in \C^N$ with the following formula(s),
\begin{equation}\label{iframeop}
    \forall  x \in \C^N, \text{ }
    x = \sum_{i=1}^M \langle x, v_i \rangle S^{-1}v_i = \sum_{i=1}^M
    \langle x, S^{-1}v_i \rangle v_i.
\end{equation}
$\mathcal{F}$ is a tight frame if and only if $S = A \,I\!d$ where $A$ is the frame bound. 
This allows for easy reconstruction of any $x \in \C^N$ and also a sufficient condition for showing
a sequence of vectors in $\C^N$ is a tight frame.


\begin{defn}
Let $\vphi \in \C^N$. The \textit{discrete short-time Fourier transform} of $\vphi$ with respect
to $\psi \in \C^N$ is defined by
\[
V_\psi(\vphi)[m,n] = \langle \vphi, e_n \tau_m \psi \rangle
= \sum_{k=0}^{N-1} \vphi[k] \overline{\psi[k-m]} e^{-2 \pi i  nk / N}
= (\vphi \tau_m(\overline{\psi}))^{\widehat{}}\,[n].
\]
The inversion formula is given by
\[
\vphi[k] = \frac{1}{N\|\psi\|_2^2} \sum_{n = 0}^{N-1}\sum_{m=0}^{N-1}
V_\psi(\vphi)[m,n] e_n \tau_m \psi.
\]
\end{defn}


The following theorem shows that full Gabor systems, or Gabor systems generated by
$\Lambda = (\Z/N\Z) \times (\Z/N\Z)^{\widehat{}}$, are always tight frames, regardless
of the choice of $\vphi$.
\begin{thm}\thlabel{fullgabor}
Let $\vphi \in \C^N \setminus \{0\}$ and let $\Lambda = (\Z/N\Z) \times (\Z/N\Z)^{\widehat{}}$. 
Then, the Gabor system $(\vphi, \Lambda)$ is a tight frame with frame bound $N\|\vphi\|_2^2$.
\end{thm}

\begin{prf}
    For every $x \in \C^N$ we compute $S(x)$,
\[
S(x) = \sum_{n = 0}^{N-1} \sum_{m = 0}^{N-1} \langle x, e_n \tau_m \vphi \rangle
    e_n \tau_m \vphi = \sum_{n=0}^{N-1}\sum_{m=0}^{N-1}V_\vphi(x)[m,n] e_n \tau_m \vphi
= N\|\vphi\|_2^2 x. \blacksquare
\]
\end{prf}

In light of \thref{fullgabor}, we only want to analyze Gabor systems where $\Lambda$ is a proper 
subgroup of 
$(\Z/N\Z) \times (\Z/N\Z)^{\widehat{}}$ since the tightness of the Gabor frame $(\vphi,\Lambda)$
is completely independent of the choice of $\vphi$ if $\Lambda = (\Z/N\Z) \times
(\Z/N\Z)^{\widehat{}}$.

The primary tool in our analysis will be Janssen's representation.
Part of Janssen's representation includes the adjoint subgroup of the subgroup $\Lambda \subseteq
(\Z/N\Z) \times (\Z/N\Z)^{\widehat{}}$, whose definition is given below.

\begin{defn}\thlabel{adjointdef}
Let $\Lambda$ be a subgroup of $(\Z/N\Z) \times (\Z/N\Z)^{\widehat{}}$. The \textit{adjoint
subgroup}, $\Lambda^\circ$, is defined by
\[
\Lambda^\circ = \{(m,n)\in(\Z/N\Z) \times (\Z/N\Z)^{\widehat{}}: e_n\tau_m e_\ell \tau_k
= e_\ell \tau_k e_n \tau_m,  \forall (k,\ell) \in \Lambda \}.
\]
\end{defn}
In other words, $\Lambda^\circ$ consists of the time-frequency operators which commute with
all time-frequency operators in $\Lambda$. 
The following form of Janssen's representation is less general than what is usually
known as Janssen's representation, but we choose to use this form because it is adjusted for use in
our main theorem. A more general version is proved in \cite{Pfa2013}.

\begin{thm}\thlabel{Janssen}
Let $\Lambda$ be a subgroup of $(\Z/N\Z) \times (\Z/N\Z)^{\widehat{}}$ and $\Lambda^\circ$ be the
adjoint subgroup of $\Lambda$. Let $\vphi \in \C^N \setminus \{0\}$. Then, the Gabor frame operator 
of the Gabor system $(\vphi, \Lambda)$ can be written as
\[
S = \frac{|\Lambda|}{N} \sum_{(m,n) \in \Lambda^\circ} \langle e_n \tau_m \vphi, \vphi 
\rangle e_n \tau_m.
\]
\end{thm}


\begin{defn}
Let $\vphi \in \C^N$. The \textit{discrete periodic ambiguity function} (DPAF) of $\vphi$ is the
function $A_p(\vphi): (\Z/N\Z) \times (\Z/N\Z)^{\widehat{}} \to \C$ defined by
\begin{equation}\label{dpafform}
A_p(\vphi)[m,n] = \frac{1}{N} \sum_{k=0}^{N-1} \vphi[k+m]\overline{\vphi[k]} e^{-2\pi i n k /N}
= \frac{1}{N} \langle \tau_{-m}\vphi, e_{n}\vphi \rangle.
\end{equation}
\end{defn}

It should be noted that the discrete periodic ambiguity function is essentially the same
as the short-time Fourier transform of $\vphi$ with $\vphi$ itself as the window function, and 
thus can be thought of as essentially interchangable. 
The computation in (\ref{dpaf2stft}) demonstrates this idea.

\begin{equation}\label{dpaf2stft}
    A_p(\vphi)[m,n] = \frac{1}{N}\langle \tau_{-m}\vphi, e_n\vphi \rangle 
    = \frac{e^{2\pi i m n / N}}{N} \langle \vphi, e_n \tau_m \vphi \rangle
    = \frac{e^{2\pi i m n / N}}{N} V_{\vphi}(\vphi)[m,n].
\end{equation}

\begin{defn}\thlabel{apsparse}
    Let $\vphi \in \C^N \setminus \{0\}$ and let $\Lambda \subseteq \Z/N\Z \times
    (\Z/N\Z)^{\widehat{}}$ be a subgroup. 
    Let $\Lambda^\circ$
    be the adjoint subgroup of $\Lambda$. $A_p(\vphi)$ is 
    \textit{$\Lambda^\circ$-sparse} if for every $(m,n) \in \Lambda^\circ \setminus \{
    (0,0) \}$ we have that $A_p(\vphi)[m,n] = 0$.
\end{defn}


We shall prove that $\Lambda^\circ$-sparsity is a necessary and sufficient condition for
determining whether or not a given Gabor system is a tight frame. To accomplish this we need 
one more theorem. 
Recall that the space of linear operators on $\C^N$ forms an $N^2$-dimensional space. 
Moreover, given any orthonormal basis $\{e_i\}_{i=1}^N$, we can define the
\textit{Hilbert-Schmidt inner product} of two linear operators $A,B$ by
\begin{equation}
    \langle A, B \rangle_{HS} = \sum_{i=1}^N \sum_{j=1}^N \langle A e_i, e_j \rangle
    \langle B e_i, e_j \rangle.
\end{equation}
The Hilbert-Schimdt inner product is independent of choice of orthonormal basis. We call the
space of linear operators on $\C^N$ equipped with the Hilbert-Schmidt inner product the
\textit{Hilbert-Schmidt space}.
\thref{linindop} is given without proof, but a proof can be found in \cite{Pfa2013}.

\begin{thm}\thlabel{linindop}
    The set of normalized time frequency translates 
    $\{ \frac{1}{\sqrt{N}}e_\ell \tau_k : (k,\ell) \in (\Z/N\Z) \times
    (\Z/N\Z)^{\widehat{}} \, \}$ forms an orthonormal basis for the $N^2$-dimensional 
    Hilbert-Schmidt space of linear operators on $\C^N$.
\end{thm}

With \thref{Janssen}, \thref{linindop}, and \thref{apsparse} we are ready to 
present the main theoretical result.


\begin{thm}\thlabel{aptightframe}
    Let $\vphi \in \C^N \setminus \{0\}$ and let $\Lambda \subseteq \Z/N\Z \times 
    (\Z/N\Z)^{\widehat{}}$. $(\vphi, \Lambda)$ is a tight frame if and only if $A_p(\vphi)$ is
    $\Lambda^\circ$-sparse. Moreover, the frame bound is given by 
    $|\Lambda| A_p(\vphi)[0,0]$.
\end{thm}

\begin{prf}
    By Janssen's representation and using the definition of $A_p(\vphi)$ we have
    \[
        S = \frac{|\Lambda|}{N} \sum_{(k,\ell) \in \Lambda^\circ} \langle e_\ell
        \tau_k \vphi, \vphi \rangle e_\ell \tau_k
        = \frac{|\Lambda|}{N} \sum_{(k,\ell) \in \Lambda^\circ} \langle 
        \tau_k \vphi, e_{-\ell} \vphi \rangle e_\ell \tau_k
    \]
    \begin{equation}\label{tightframejan}    
        = |\Lambda| \sum_{(k,\ell) \in \Lambda^\circ} A_p(\vphi)[-k,-\ell] e_\ell \tau_k
        = |\Lambda| \sum_{(k,\ell) \in \Lambda^\circ} A_p(\vphi)[k,\ell] e_\ell \tau_k
    \end{equation}
    where the last equality comes from the fact that $\Lambda^\circ$ is a subgroup of
    $(\Z/N\Z) \times (\Z/N\Z)^{\widehat{}}$. Clearly,
    if $A_p(\vphi)$ is $\Lambda^\circ$-sparse, then by (\ref{tightframejan}) the frame
    operator will be $|\Lambda| A_p(\vphi)[0,0]$ times the identity. It remains to show
    that $A_p(\vphi)$ is a necessary condition. For $S$ to be tight we need 
    \begin{equation}\label{jansencond}
        S = |\Lambda| \sum_{(k,\ell) \in \Lambda^\circ} A_p(\vphi)[k,\ell] e_\ell \tau_k
        = A \, I\!d.
    \end{equation}
    In particular, we can rewrite (\ref{jansencond}) to
    \[
        \sum_{(k,\ell) \in \Lambda^\circ \setminus \{(0,0)\}} |\Lambda|A_p(\vphi)[k,\ell] e_\ell 
        \tau_k
        + (|\Lambda|A_p(\vphi)[0,0] - A) I\!d = 0.
    \]
    By \thref{linindop}, the set of linear operators $\{e_\ell \tau_k 
    : (k,\ell) \in \Lambda^\circ \} $ is linearly independent. Thus, $A_p(\vphi)[k,\ell] = 0$ for
    every $(k,\ell) \in \Lambda^\circ \setminus \{(0,0)\}$ and $A = |\Lambda|A_p(\vphi)[0,0]$.
    We conclude that, $A_p(\vphi)$ is $\Lambda^\circ$-sparse and the frame has the desired
    frame bound. $\blacksquare$
\end{prf}


\thref{aptightframe} is closely connected to the Wexler-Raz criterion. The Wexler-Raz criterion
checks whether a Gabor system $(\tilde{\vphi},\Lambda)$ is a dual frame
to $(\vphi,\Lambda)$. In particular, if $S$ is the frame operator of $(\vphi,\Lambda)$, then
$S^{-1}\vphi$ is the canonical dual frame of $(\vphi,\Lambda)$. 
\thref{aptightframe} is a special case of Wexler-Raz which confirms that the
canonical dual frame associated with tight frames indeed satisfies the Wexler-Raz criterion.
Again, a proof is not given but can be found in \cite{Pfa2013}.

\begin{thm}\thlabel{wexraz}
    Let $\Lambda$ be a subgroup of $(\Z/N\Z) \times (\Z/N\Z)^{\widehat{}}$. For Gabor systems
    $(\vphi,\Lambda)$ and $(\tilde{\vphi},\Lambda)$ we have
    \begin{equation}\label{dualframe}
        x = \sum_{(k,\ell) \in \Lambda} \langle x, e_\ell \tau_k \tilde{\vphi} \rangle
        e_\ell \tau_k \vphi, \hs x \in \C^N,
    \end{equation}
    if and only if
    \begin{equation}\label{wexrazcrit}
        \langle \vphi, e_\ell \tau_k \tilde{\vphi} \rangle = 
        \frac{N}{|\Lambda|} \delta_{(k,\ell),(0,0)},
        \hs (k,\ell) \in \Lambda^\circ.
    \end{equation}
\end{thm}

If $(\vphi, \Lambda)$ is a tight frame, then the cannonical dual frame to $(\vphi, \Lambda)$ is
$(A^{-1}\vphi,\Lambda)$, since $S^{-1} = A^{-1} I\!d$. If we use the Wexler-Raz criteron to
verify if $(A^{-1}\vphi,\Lambda)$ is a dual frame and use $A = |\Lambda|A_p(\vphi)[0,0]$ as
in \thref{aptightframe}, then (\ref{dualframe}) holds if and only if,

\[
    \frac{1}{|\Lambda| A_p(\vphi)[0,0]} \langle \vphi, e_\ell \tau_k \vphi \rangle
    = \frac{N}{|\Lambda|} = \delta_{(k,\ell),(0,0)}, \hs (k,\ell) \in \Lambda^\circ.
\]
This can be rewritten as
\[
    A_p(\vphi)[k,\ell] = A_p(\vphi)[0,0]\delta_{(k,\ell),(0,0)} \hs (k,\ell) \in \Lambda^\circ.
\]
This is the same condition as $A_p(\vphi)$ being $\Lambda^\circ$-sparse.

We close this section with a few operations on $\vphi$ which preservere the
$\Lambda^\circ$-sparsity of the ambiguity function. 

\begin{prop}\thlabel{apsparseprop}
    Let $\vphi \in \C^N \setminus \{0\}$ and let $\Lambda \subseteq \Z/N\Z \times
    (\Z/N\Z)^{\widehat{}}$ be a subgroup. 
    Suppose that 
    $A_p(\vphi)$ is $\Lambda^\circ$-sparse. Then the following are also $\Lambda^\circ$-
    sparse:
    \begin{itemize}
        \item[(i)] $\forall k \in \Z/N\Z$, $A_p(\tau_k\vphi)$
        \item[(ii)] $\forall \ell \in (\Z/N\Z)^{\widehat{}}$, $A_p(e_\ell \vphi)$
        \item[(iii)] $\forall c \in \C\setminus\{0\}$, $A_p(c \vphi)$
        \item[(iv)] $A_p(\overline{\vphi})$.
    \end{itemize}
\end{prop}

\begin{prf}
    Each part follows from direct computation which is shown below:

    \begin{itemize}
        \item[(i)] Let $k \in (\Z/N\Z)$. Then,
        \[
            A_p(\tau_k \vphi)[m,n] = \frac{1}{N} \sum_{j = 0}^{N-1} \vphi[j-k+m] \overline{
            \vphi[j-k]} e^{-2\pi i n j / N}
        \]
        \[
            = \frac{1}{N}\sum_{j' = -k}^{N-k-1} \vphi[j'+m] \overline{\vphi[j']} e^{-2 \pi i n
            (j' + k) / N} 
            = e^{-2 \pi i n k / N} 
            \frac{1}{N}\sum_{j' = 0}^{N-1} \vphi[j'+m] \overline{\vphi[j']} e^{-2 \pi i n
            j'  / N} 
        \]
        \[
            = e^{-2 \pi i n k / N} A_p(\vphi)[m,n] = 0
        \]
        if $A_p(\vphi)[m,n] = 0$.
        
        \item[(ii)] Let $\ell \in (\Z/N\Z)^{\widehat{}}$. Then,
        \[
            A_p(e_\ell \vphi)[m,n] = \frac{1}{N} \sum_{j = 0}^{N-1} e^{2 \pi i \ell(j+m)/ N}
            \vphi[j+m]
            e^{-2 \pi i \ell j / N} \overline{\vphi[j]}
            e^{-2 \pi i n j / N}
        \]
        \[
            = e^{2 \pi i \ell m / N} \frac{1}{N} \sum_{j = 0}^{N-1} \vphi[j+m]
            \overline{\vphi[j]} e^{-2 \pi i n j / N}
            = e^{2 \pi i \ell m / N} A_p(\vphi)[m,n] = 0
        \]
        if $A_p(\vphi)[m,n] = 0$.
        \item[(iii)] Let $c \in \C\setminus\{0\}$. Then,
        \[
            A_p(c\vphi)[m,n] = \frac{1}{N} \sum_{j=0}^{N-1} c\vphi[j+m] \overline{c}
            \overline{\vphi[j]} e^{-2 \pi i n j / N}a = |c|^2 A_p(\vphi)[m,n] = 0
        \]
        if $A_p(\vphi)[m,n] = 0$.
        \item[(iv)] By direct computation,
        \[
            A_p(\overline{\vphi})[m,n] = \frac{1}{N} \sum_{j=0}^{N-1} \overline{\vphi[j+m]}
            \vphi[j] e^{-2 \pi i n j / N} = \frac{1}{N} \sum_{j' = m}^{N + m - 1}
            \vphi[j'] \overline{\vphi[j'+m]} e^{-2 \pi i n (j' - m) / N} 
        \]
        \[
            = e^{2 \pi i n m /N} \frac{1}{N} \sum_{j' = 0}^{N-1} \vphi[j'+m]
            \overline{\vphi[j']} e^{-2 \pi i n j' / N} 
        \]
        \[
            = e^{2 \pi i n m /N} A_p(\vphi)[m,n] = 0
        \]
        if $A_p(\vphi)[m,n] = 0$. $\blacksquare$
    \end{itemize}
\end{prf}


\section{CAZAC Sequences}\label{cazacsec}
Since the goal is to analyze which CAZACs are suitable for genetating tight Gabor frames,
we briefly discuss the background of CAZAC sequences.
A more detailed exposition on CAZACs can be found in \cite{BenCorMag2016}.
CAZAC is an acronym which stands for
Constant Amplitude and Zero Autocorrelation. These sequences have applications in areas
such as coding theory \cite{LevMoz2004} and have several interesting mathematical properties
as well as problems. 
First, we begin with the definition of a CAZAC sequence.


\begin{defn}\thlabel{cazacdef}
    $\vphi \in \C^N$ is a \textit{CAZAC sequence} if the following two properties hold:
\begin{equation}\tag{CA}
    |\vphi[k]| = 1, \hs \forall k \in (\Z/N\Z)
\end{equation}
and
\begin{equation}\tag{ZAC}
    \sum_{k=0}^{N-1} \vphi[k+m] \overline{\vphi[k]} = 0, \hs \forall m \in (\Z/N\Z) \setminus
    \{0\}.
\end{equation}
\end{defn}

An equivalent definition of CAZAC sequences are sequences where
both $\vphi$ and $\widehat{\vphi}$ satisfiy (CA). 
As such, CAZAC sequences are sometiems refered to biunimodular sequences.
This idea is made clearer by \thref{cazacprop}.
There are seven CAZAC sequences which shall be used in the examples: Chu, P4, Wiener, square-length
Bj\"{o}rck-Saffari, Milewski, and Bjorck sequences. We shall define these sequences now and leave 
the analysis of their DPAFs to Section \ref{cazacdpaf}.


\begin{defn}
Let $\vphi \in \C^N$ be of the form
\[
\vphi[k] = e^{\pi i p[k] / N},
\]
where $p[k]$ is a polynomial. Then we can define the Chu, P4, and Wiener CAZAC sequences with,
\begin{itemize}
    \item Chu: $p[k] = k (k-1)$ if $N$ is odd.
    \item P4: $p[k] = k (k-N)$ for any $N$.
    \item Wiener (odd length): $p[k] = sk^2$, if $N$ is odd and gcd$(s,N) = 1$.
    \item Wiener (even length): $p[k] = sk^2/2$, if $N$ is even and gcd$(s,2N) = 1$.
\end{itemize}
\end{defn}

All three of these sequences belong to a class of sequences known as chirp sequences. Chirp
sequences are sequence whose frequency change linearly in time. They are also known as
quadratic phase sequences. Some additional exposition about these sequences can be found in
\cite{BenKonRan2009}.

\begin{defn}
    Let $c \in \C^N$ be unimodular, $v \in \C^M$ be CAZAC, and $\sigma$ be any permutation of the
    set $\{0,1,\cdots,N-1\}$. Then,
    \begin{itemize} 
        \item[(i)] The \textit{square length Bj\"{o}rck-Saffari CAZAC sequence},
        $\vphi \in \C^{N^2}$, is defined by,
        \[
            \forall r \in (\Z/N\Z), h \in (\Z/N\Z), \vphi[aN+b] = c[h]e^{2\pi i r \sigma(h) /N}.
        \]
    
    \item[(ii)] The \textit{Milewski CAZAC sequence}, $\vphi \in \C^{MN^2}$ is defined by,
        \[
            \forall a \in (\Z/MN\Z), b \in (\Z/N\Z), \vphi[aN+b] = v[a] e^{2 \pi i a b / MN}.
        \]
    \end{itemize}
\end{defn}

The square-length Bj\"orck-Saffari sequences were defined by Bj\"orck and Saffari as a building
block to a more general class of CAZAC sequences whose length were not necessarily a perfect
square \cite{BjoSaf1995}. The Milewski sequence was first defined by Milewski in \cite{Mil1983}
as a way to construct more CAZAC sequences out of already existing ones. 


\begin{defn}
    Let $p$ be prime. The \textit{Legendre symbol} is defined by
    \[
    \left( \frac{k}{p} \right) =
    \begin{cases}
    0, & \text{if }k \equiv 0 \mod p \\
    1, & \text{if }k \equiv n^2 \mod p \text{ for some } n \in (\Z/p\Z) \setminus \{0\} \\
    -1, & \text{otherwise.}
    \end{cases}
    \]
\end{defn}

\begin{defn}\thlabel{bjodef}
    Let $p$ be prime and let $\vphi \in \C^p$ be of the following form,
    \begin{equation}\label{bjoform}
    \vphi[k] = e^{i \theta[k]}.
    \end{equation}
    We define the \textit{Bj\"orck sequence} by letting $\theta[k]$ in (\ref{bjoform}) be as
    follows,
    \begin{itemize}
    
        \item If $p \equiv 1 \mod 4$ then,
            \[
                \theta[k]= \left(\frac{k}{p}\right) \text{arccos}\left(\frac{1}
                {1+\sqrt{p}}\right).
            \]
        \item If $p \equiv -1 \mod 4$ then,
            \[
            \theta[k] = 
            \begin{cases}
                \text{arccos}\left(\frac{1-p}{1+p}\right), & \text{if } 
                \left(\frac{k}{p}\right) = -1 \\
                0 & \text{otherwise.}
            \end{cases}
            \]
    \end{itemize}
\end{defn}

\begin{figure}[h!]
\centering
\includegraphics[width=0.9\textwidth]{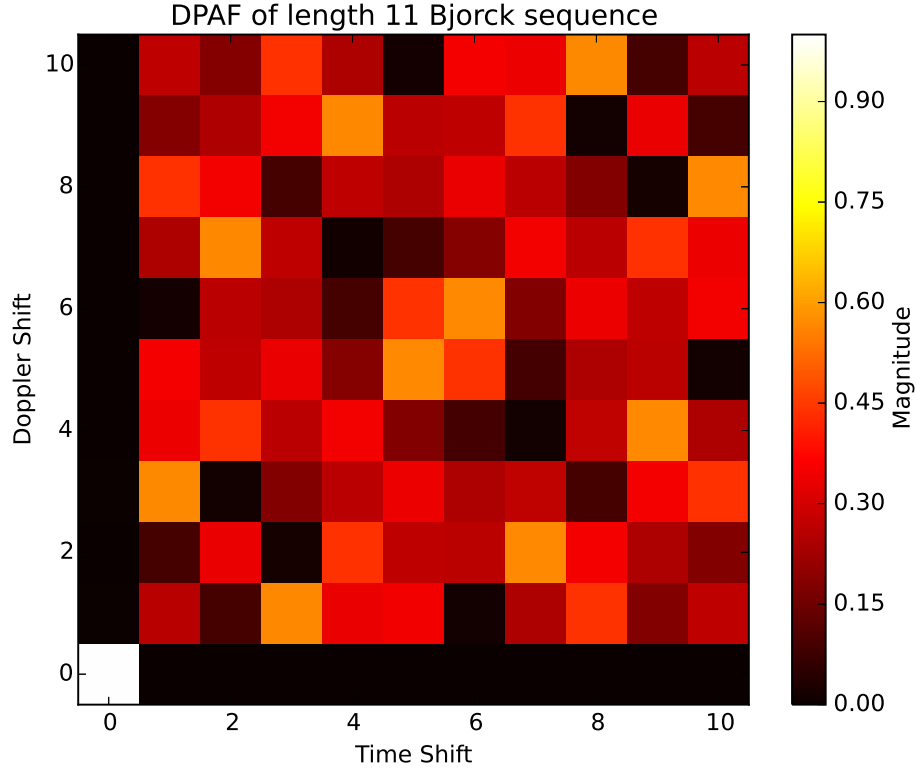}
\caption{DPAF of length 11 Bj\"orck sequence.}
\label{bjo11}
\end{figure}

Figure \ref{bjo11} shows that the length 11 Bj\"orck sequence is 
indeed a CAZAC sequence. In fact, all sequences generated by \thref{bjodef} are CAZAC 
\cite{Bjo1990}. 
In this particular example, one can see that the discrete periodic ambiguity function of the
length 11 Bj\"orck sequence is almost always nonzero,
except for on $A_p(\vphi)[\cdot,0]$ and $A_p(\vphi)[0,\cdot]$. 
These properties hold true for any Bj\"{o}rck sequence of any prime length.
In light of \thref{aptightframe}, we can see that Bj\"{o}rck sequences are ill-suited to the
construction of tight frames, despite being CAZAC. We shall explore this idea more in
Section \ref{exsparse}.
More detailed exploration on the behavior of the DPAF of the Bj\"orck sequence can be found
in \cite{BenBenWoo2012} and \cite{KebKonBen2007}.
For completeness, the length 11 Bj\"orck
sequence is listed out below,
\[
    \vphi = (1, 1, e^{i\theta_{11}}, 1, 1, 1, e^{i\theta_{11}}, e^{i\theta_{11}}, e^{i\theta_{11}}, 
    1, e^{i\theta_{11}}), 
\]
where $\theta_{11} = \text{arccos}(-10/11)$.


\begin{prop}\thlabel{cazacprop}
    Let $\vphi \in \C^N$ and let $c \in \C$ be such that $|c| = 1$. Then,
    \begin{itemize}
        \item[(i)] If $\vphi$ is CA, then $\widehat{\vphi}$ is ZAC.
        \item[(ii)] If $\vphi$ is CAZAC, then $\widehat{\vphi}$ is also CAZAC.
        \item[(iii)] If $\vphi$ is CA, then $\widehat{\vphi}$ can have zeros.
        \item[(iv)] If $\vphi$ is CAZAC, then $c\vphi$ is also CAZAC.
    \end{itemize}
\end{prop}

Property (iv) allows us to say two CAZACs are equivalent if they are complex rotations of each
other. We can assume that the representative CAZAC in each equivalence class is the sequence whose 
first entry is 1. Given this, a natural question is as follows: For each $N$, how many CAZACs are 
there of length $N$? If $N$ is a prime number, then there
are only finitely many classes of CAZAC sequences \cite{Haa2008}. On the other
hand, if $N$ is composite and is not square-free, then there are infinitely many classes
of CAZAC sequences \cite{BjoSaf1995}. If $N$ is composite and square-free, it is unknown
whether the number of equivalence classes is finite or infinite. 
Another question is, what other CAZACs can we construct besides the Chu, P4, etc.?
The discovery of other CAZAC sequences can be transformed into two different
problems in very different areas of mathematics. We shall explore
these two other equivalent problems: The first involves circulant Hadamard matrices and the
second is with cyclic n-roots.


\begin{defn}\thlabel{hadacircdef}
    Let $H \in \C^{N \times N}$. 
    \begin{itemize}
        \item[(i)] $H$ is a \textit{Hadamard matrix} if for every $(i,j)$, $|H_{i,j}| = 1$
        and $H H^* = N I\!d$.
        \item[(ii)] $H$ is a \textit{circulant matrix} if for each $i$, the $i$-th row is
        a circular shift of the first row by $i-1$ entries to the right.
    \end{itemize}
\end{defn}

One can construct a circulant Hadamard matrix by making the first row a CAZAC sequence and
each row after a shift of the previous row to the right. It is clear that this construction
leads to a circulant matrix, and the ZAC property will guarantee that $H H^* = N I\!d$.
Thus, given a CAZAC, we can generate a Hadamard matrix, and even better, we have the following
theorem \cite{BenKonRan2009}.

\begin{thm}\thlabel{hadacirc}
    $\vphi \in \C^{N \times N}$ if and only if the circulant matrix generated by $\vphi$ is
    a Hadamard matrix.
\end{thm}

In particular, \thref{hadacirc} gives a one-to-one correspondence between circulant $N \times N$ 
Hadamard matrices whose diagonal consists of ones and the equivalence classes of CAZAC sequences 
of length $N$. Therefore, an equiavlent problem to the discovery of additional CAZAC sequences is
the computation of circulant Hadamard matrices. 
There is a significant amount of research and interest in Hadamard matrices,
even outside of the context of CAZACs. A catalogue of
complex Hadamard matrices with relevant citations can be found online at \cite{BruTadWeb}.


\begin{defn}
$(z_0, z_1, \cdots, z_{N-1}) \in \C^N$ is a \textit{cylclic N-root} if it satisfies the following
system of equations,
\begin{equation}\label{cyclicnroot}
\begin{cases}
z_0 + z_1 + \cdots + z_{N-1} = 0 \\
z_0z_1 + z_1z_2 + \cdots + z_{N-1}z_0 = 0 \\
z_0z_1z_2 + z_1z_2z_3 + \cdots + z_{N-1}z_0z_1 = 0 \\
\cdots \cdots \\
z_0z_1z_2\cdots z_{N-1} = 1
\end{cases}.
\end{equation}
\end{defn}

Let $\vphi \in \C^N$ be CAZAC with $\vphi[0] = 1$. Let us emphasize the sequence nature of
$\vphi \in \C^N$ and write $\vphi[k] = \vphi_k$. Then, 
\begin{equation}\label{cazac2croot}
(z_0, z_1, \cdots, z_{N-1}) :=
\left( \frac{\vphi_1}{\vphi_0}, \frac{\vphi_2}{\vphi_1}, \cdots, 
\frac{\vphi_0}{\vphi_{N-1}} \right)
\end{equation}
is a cylcic n-root. In fact, there is a one-to-one correspondence between cyclic n-roots
and CAZAC sequences whose first entry is one \cite{Haa2008}. In the same manner as circulant
Hadamard matrices, finding cyclic N-roots is equivalent to finding CAZAC sequences of length $N$.
In particular, using (\ref{cazac2croot}) we can see that we can construct CAZAC sequences by
the following (recursive) formula:
\begin{equation}\label{croot2cazac}
\forall k \in \{1,2,\cdots,N-1\}, \text{ } \vphi_0 = 1, \text{ } \vphi_k = \vphi_{k-1} z_{k-1}.
\end{equation}
Cyclic $n$-roots can be used to show that the number of prime length CAZACs is finite. This
was proved by Haagerup \cite{Haa2008}, but a summary, along with more results on CAZACs, can be
found in \cite{BenCorMag2016}.

It is still unclear what the exact role of CAZACs is in the generation of tight frames. 
Many of the known CAZAC sequences are suitable for generating tight Gabor frames, as we shall
see in Section \ref{exsparse}. On the other hand, the Bj\"{o}rck sequence is also CAZAC but is
very ill-suited for generating tight Gabor frames.
Part of the difficulty is that although there are results quantifying the number of
CAZAC sequences for given length $N$, very few of them have been explicitly written out.
Many of the ones which are known are generated by roots of unity. The Bj\"{o}rck sequence is
the exception to this and is also the one that is ill-suited for generating tight Gabor frames.
One could perhaps show that all CAZAC sequences generated by roots of unity
(eg. Chu, P4, Wiener, roots of unity generated Milewski) will have sparse discrete periodic
ambiguity functions. Hopefully, the eventual discovery of more CAZAC sequences will help
further clarify the usability of CAZAC sequences to generate tight Gabor frames.

\section{Discrete periodic ambiguity functions of selected sequences}\label{cazacdpaf}

The following two sections will be devoted to examples of \thref{aptightframe} using
the CAZAC sequences from Section \ref{cazacsec} and various time-frequency subgroups. 
In this section we will compute the DPAFs for five classes of sequences: Chu, P4, Wiener,
Square-length Bj\"orck-Saffari, and Milewski sequences. 

\subsection{Chu Sequence} \label{chudpaf}
The computation of the DPAF of the Chu sequence is as follows,

\[
A_p(\vphi)[m,n] = \frac{1}{N} \sum_{k=0}^{N-1} e^{\pi i (k+m)(k+m-1)/N} e^{-\pi i k(k-1)/N} 
e^{-2 \pi i n k / N}
= \frac{1}{N} \sum_{k=0}^{N-1} e^{\pi i [(k+m)(k+m-1)-k(k-1)-2nk] / N}
\]
\[
= \frac{1}{N} \sum_{k=0}^{N-1} e^{\pi i (k^2 + 2km + m^2 - k - m - k^2 + k - 2nk) / N}
= \frac{1}{N} \sum_{k=0}^{N-1} e^{\pi i (2km + m^2 - m  - 2nk) / N}
\]
\[
= \frac{1}{N} e^{\pi i (m^2-m) / N} \sum_{k=0}^{N-1} e^{2 \pi i k (m-n)/N}
\]
\[
= 
\begin{cases}
e^{\pi i (m^2-m)/N}, & \text{if } m \equiv n \mod N \\
0, & \text{otherwise.}
\end{cases}
\]
Note that in particular, everything off of the line $n=m$ returns a zero. We will leverage
this fact in several of the examples in Section \ref{exsparse}.

\subsection{P4 Sequence} \label{p4dpaf}
The computation of the DPAF of the P4 sequence is as follows,
\[
    A_p(\vphi)[m,n] = \frac{1}{N} \sum_{k=0}^{N-1} e^{\pi i (k+m)(k+m-N) / N}
    e^{-\pi i k(k-N) / N} e^{-2 \pi i n k /N}
\]
\[
    = \frac{1}{N} \sum_{k=0}^{N-1} e^{\pi i [(k+m)(k+m-N) - k(k-N) - 2nk] / N}
    = \frac{1}{N} \sum_{k=0}^{N-1} e^{\pi i (k^2 + 2km + m^2 - Nk - Nm - k^2 + Nk - 2nk) / N}
\]
\[
    = \frac{1}{N} \sum_{k=0}^{N-1} e^{\pi i (m^2 - Nm + 2mk - 2nk) / N}
    = \frac{1}{N}(-1)^m e^{\pi i m^2 / N} \sum_{k=0}^{N-1} e^{2 \pi i (m-n) k / N}
\]
\[
=
\begin{cases}
    (-1)^m e^{\pi i m^2 / N}, & \text{if } m \equiv n \mod N \\
    0, & \text{otherwise.}
\end{cases}
\]
Like the Chu sequence, the DPAF of the P4 sequence is also only nonzero on the diagonal
$n = m$. Due to this fact, the Chu and P4 sequences will be used interchangably in the examples.

\subsection{Wiener Sequences}
We start with the case where $N$ is odd. In this case, $\vphi$ has the form
\[
    \forall k \in (\Z/N\Z), \text{ } \vphi[k] = e^{2 \pi i s k^2 / N}
\]
Then, the computation of the DPAF is as follows,
\[
    A_p(\vphi)[m,n] = \frac{1}{N} \sum_{k=0}^{N-1} e^{2 \pi i s (k+m)^2 / N}
    e^{-2 \pi i s k^2 / N} e^{-2 \pi i n k / N}
    = \frac{1}{N} \sum_{k=0}^{N-1} e^{2 \pi i (sk^2 + 2skm + sm^2 - sk^2 - nk) / N}
\]
\[
    = \frac{1}{N} \sum_{k=0}^{N-1} e^{2 \pi i (sm^2 + 2skm - nk)/N}
    = \frac{1}{N} e^{2 \pi i s m^2 / N} \sum_{k=0}^{N-1} e^{2\pi i (2sm-n) k /N}
\]
\[
=
\begin{cases}
    e^{2 \pi i s m^2 / N}, & \text{if } 2sm \equiv n \mod N \\
    0, & \text{otherwise.}
\end{cases}
\]
The second case is where $N$ is even. In this case, $\vphi$ has the form
\[
    \vphi[k] = e^{\pi i s k^2 / N}, \hs[3mm] k \in \{0, \cdots, N-1\}.
\]
In this case, the computation of the DPAF is as follows,
\[
    A_p(\vphi)[m,n] = \frac{1}{N} \sum_{k=0}^{N-1} e^{\pi i s (k+m)^2 / N} e^{-\pi i s k^2 / N}
    e^{-2 \pi i n k / N} = 
    \frac{1}{N} \sum_{k=0}^{N-1} e^{\pi i (sk^2 + 2skm + sm^2 - sk^2 - 2nk) / N}
\]
\[
    = \frac{1}{N} \sum_{k=0}^{N-1} e^{\pi i (sm^2 + 2skm - 2kn) / N} = e^{\pi i sm^2 / N} 
    \frac{1}{N} \sum_{k=0}^{N-1} e^{2 \pi i (sm-n) k /N}
\]
\[
=
\begin{cases}
e^{\pi i s m^2 / N}, & \text{if } sm \equiv n \mod N \\
0, &  \text{otherwise.}
\end{cases}
\]

\subsection{Square Length Bj\"orck-Saffari Sequences} \label{slbsdpaf}
The computation for the DPAF of the square length Bj\"orck-Saffari sequences is as follows,
\[
A_p(\vphi)[sN+t,kN+\ell]=\frac{1}{N^2} \sum_{r=0}^{N-1}\sum_{h=0}^{N-1}
\vphi[(r+s)N+(h+t)]\overline{\vphi[rN+h]}
e^{-2 \pi i (kN+\ell)(rN+h)/N^2}
\]
\[
= \frac{1}{N^2} \sum_{r=0}^{N-1}\sum_{h=0}^{N-1}
c[h+t]\overline{c[h]}e^{2\pi i (r+s + \floor{\frac{(h+t)}{N}})\sigma(h+t)/N}
e^{-2 \pi i r \sigma(h) / N} e^{-2 \pi i [(kh+r\ell)N+\ell h]/N^2}
\]
\[
= \frac{1}{N^2} \sum_{h=0}^{N-1}
c[h+t]\overline{c[h]}e^{2\pi i (s + \floor{\frac{(h+t)}{N}})\sigma(h+t)/N}
e^{-2 \pi i (khN+\ell h)/N^2}
\sum_{r=0}^{N-1}e^{2 \pi i (\sigma(h+t)-\sigma(h)-\ell)r/N}
\]
\[
= \frac{1}{N} \sum_{h=0}^{N-1}
c[h+t]\overline{c[h]}e^{2\pi i (s + \floor{\frac{(h+t)}{N}})\sigma(h+t)/N}
e^{-2 \pi i (khN+\ell h)/N^2} 
\]
if $\sigma(h+t) - \sigma(h) - \ell \equiv 0 \mod N$, and 0 otherwise.
In particular, if $\sigma(h) = h$ for all $h$, then the above condition reduces to $t \equiv \ell
\mod N$.


\subsection{Milewski sequences}\label{mildpaf}
We shall write out the DPAF of the Milewski sequence without computation, but the computation 
can be found in \cite{BenDon2007}. The DPAF of the Milewski sequence is,
\begin{equation}\label{milaf}
    N A_p(\vphi)[kN+\ell,s] = \begin{cases}
        0, & \text{ if } m \not\equiv n \mod N \\
        \sum_{j=0}^{N-1} e^{2 \pi i \left(k + \lfloor\frac{j+\ell}{N}\rfloor - js\right)/MN}
        A_p(v)\left[ k + \lfloor \frac{j+\ell}{N} \rfloor, 
        \frac{s - \ell}{N} \right], & \text{ if } m \equiv  n \mod N. \\ 
    \end{cases}
\end{equation}

\section{Examples of tight Gabor frames generated by $\Lambda^\circ$-sparsity} \label{exsparse}

For the first examples, we will use the following subgroup $\Lambda \subseteq
(\Z/N\Z) \times (\Z/N\Z)^{\widehat{}}$: Let $K = \{0,a,\cdots,(bN'-1)a\}$ and $L = \{0,b,\cdots,
(aN'-1)b\}$ where $N = abN'$ and gcd$(a,b) = 1$ and let $\Lambda = K \times L$. 
We first shall compute the adjoint subgroup of $\Lambda = K \times L$.
This requires us to compute which time-frequency translates $(m,n)$ commute with every 
time-frequency translate in $K \times L$. To that end, let $(k,\ell)$ and $(m,n)$
be two time-frequency translates. We compute,
\[
(e_n \tau_m e_\ell \tau_k) \vphi[j] = (e_n \tau_m) (e^{2 \pi i \ell j / N} \vphi[j-k])
= e^{2 \pi i n j /N} e^{2 \pi i \ell (j-m) / N} \vphi[j-k-m]
\]
\[
= e^{-2 \pi i \ell m / N} (e_{\ell+n} \tau_{k+m}) \vphi[j]
\]
and
\[
(e_\ell \tau_k e_n \tau_m) \vphi[j] = (e_\ell \tau_k) (e^{2 \pi i n j / N} \vphi[j-m])
= e^{2 \pi i \ell j /N} e^{2 \pi i n (j-k) / N} \vphi[j-k-m]
\]
\[
= e^{-2 \pi i n k /N } (e_{\ell+n} \tau_{k+m}) \vphi[j].
\]
Thus, $(m,n) \in \Lambda^\circ$ if and only if
\begin{equation}\label{kladjoint}
    \ell m \equiv kn \mod N, \hs \forall (k,\ell) \in \Lambda.
\end{equation}
Since $k \in K$ and $\ell \in L$ we can write $k = k'a$ and $\ell = \ell'b$ for some $k' \in \{0,
\cdots,bN'-1\}$ and $\ell' \in \{0,\cdots,aN'-1\}$. Using this in (\ref{kladjoint}) gives a new
condition
\begin{equation}\label{kladjoint2}
\ell' b m \equiv k' a n \mod N, \forall (k',\ell') \in \{0,\cdots,bN'-1\} \times
\{0,\cdots,aN'-1\}.
\end{equation}


\begin{lemma}\thlabel{72sol}
Let $N = abN'$ where gcd$(a,b) = 1$.
$(m,n)$ is a solution to (\ref{kladjoint2}) if and only if $m$ is a multiple of $N'a$ and $n$ is
a multiple of $N'b$. In particular, if $\Lambda = K \times L$ where $K = \{0,a,\cdots,(bN'-1)a\}$
and $L = \{0,b,\cdots,(aN'-1)b\}$, then $\Lambda^\circ 
= N'K \times N'L = \{0,N'a,\cdots,(b-1)N'a\}
\times \{0,N'b,\cdots,(a-1)N'b\}$.
\end{lemma}

\begin{prf}
First note that if $m = rN'a$ and $n = sN'b$ for some $r \in \{0,\cdots,b-1\}$ and $s \in \{0, \cdots,
a-1\}$ then the left hand side of (\ref{kladjoint2}) becomes
\begin{equation}\label{72lhs}
\ell' b (rN'a) \equiv \ell' r (abN') \equiv 0 \mod N.
\end{equation}
The right hand side of (\ref{kladjoint2}) becomes
\begin{equation}\label{72rhs}
k' a (sN'b) \equiv k' s (abN') \equiv 0 \mod N.
\end{equation}
Since (\ref{72lhs}) and (\ref{72rhs}) are equal we have that 
$N'K \times N'L \subseteq \Lambda^\circ$. To show the converse, first note that since 
gcd$(a,b) = 1$, $m$ must be a multiple of $a$ and $n$ must be a multiple of $b$. In other words,
$\Lambda^\circ \subseteq K \times L$.
Now, suppose $m = ra$ and $n = sb$.
Then, condition (\ref{kladjoint2}) becomes
\begin{equation}\label{72adj}
\ell' r a b  \equiv k' s ab \mod N, \forall (k',\ell') \in \{0,\cdots,bN'-1\} \times
\{0,\cdots,aN'-1\}.
\end{equation}
If $N' = 1$, the above is always true since $N = abN' = ab$ and thus $\Lambda^\circ 
= N'K \times N'L
= K \times L$. Assume that $N' > 1$ and without loss of generality, assume $N' \nmid r$. Then,
choose $\ell' = 1$ and $k' = 0$ and it is clear that the right hand side is 0 while the left hand
side cannot be a multiple of $N$ and so (\ref{72adj}) cannot hold. Thus, $N' \mid r$ and
$N' \mid s$ is necessary and we now get $\Lambda^\circ \subseteq N'K \times N'L$ which completes
the proof. $\blacksquare$
\end{prf}

For our first example, we will use the $K \times L$ setup as above and apply it to the Chu and 
P4 sequence. As seen in Section \ref{chudpaf} and Section \ref{p4dpaf}, $A_p(\vphi)[m,n]$ only has
nonzero entries along the diagonal $m = n$ and we will leverage this fact for an easy proof of
\thref{chup4ap}. It should also be noted that some slight modifications can be made to easily
extend \thref{chup4ap} to Wiener sequences.

\begin{prop}\thlabel{chup4ap}
    Let $\vphi \in \C^N$ be either the Chu or P4 sequence and let 
    $K = \{0,a,\cdots,(bN'-1)\}$ and $L = \{0,b,\cdots,(aN'-1)b\}$ 
    with gcd$(a,b) = 1$ and $N = abN'$.  Then, the Gabor system $(\vphi, K \times L)$ is a 
    tight Gabor frame with frame
bound $NN'$.
\end{prop}

\begin{prf}
    Using \thref{72sol} we have that $(m,n) \in (K \times L)^\circ$ if and only if $m = r a N'$
    and $n = s b N'$ for some $r \in \{0,1,\cdots,bN'-1\}$ and $s \in \{0,1,\cdots,aN'-1\}$.
    However, $A_p(\vphi)[m,n] \neq 0$ if and only if $m \equiv n \mod N$. Thus, 
    $A_p(\vphi)[m,n] \neq 0$ if and only if $m,n \in N'K \cap N'L$. This intersection is
    generated by lcm$(N'a,N'b) = N'$ lcm$(a,b) = abN' = N$. From this we conclude that for
    $(m,n) \in (K \times L)^\circ$,
    $A_p(\vphi)[m,n] \neq 0$ if and only if $m = n = 0$. Thus by \thref{aptightframe}, we have
    that $(\vphi, K\times L)$ is a tight frame with frame bound $|\Lambda| = NN'$. $\blacksquare$
\end{prf}

    It should be noted that other subgroups $\Lambda \subseteq (\Z/N\Z) \times 
    (\Z/N\Z)^{\widehat{}}$ can be used in \thref{chup4ap} to get
    the same result. In this case, it is only necessary that for every 
    $(m,n) \in \Lambda^\circ$, $m \not\equiv n \mod N$ unless $(m,n) = (0,0)$.
    We demonstrate this idea with \thref{chugenframe}.

    \begin{prop}\thlabel{chugenframe}
    Let $\vphi \in \C^N$ be either the Chu or P4 sequence and let 
    $\Lambda = \{(0,0),(a,b), \\
    (2a,2b),\cdots,((N-1)a,(N-1)b)\}$ where $a \neq b$ and gcd$(a,N) = $ gcd$(b,N) = 1$.
    Then, the Gabor system $(\vphi,\Lambda)$ is a tight Gabor frame with frame bound
    $N$.
    \end{prop}

    \begin{prf}
    Suppose gcd$(a,b) = d$. Let $a = da'$ and $b = db'$. Since gcd$(a,N) = 1$ and gcd$(b,N)=1$,
    the order of the cyclic subgroup $\Lambda$ is $N$ and in particular is the same subgroup
    as the one generated by $(a',b')$. Thus, without loss of generality, we can assume that
    gcd$(a,b) = 1$. Then, by (\ref{kladjoint}), we know that $(m,n) \in \Lambda^\circ$ if and
    only if
    \[
        k n \equiv \ell m \mod N
    \]
    which can be rewritten as
    \[
        s a n \equiv  s b m \mod N.
    \]
    Letting $s = 1$, we see that $an \equiv bm \mod N$ is a necessary condition and since
    gcd$(a,b) = 1$ we need that $n$ is a multiple of $b$ and $m$ is a multiple of $a$. It is
    clear that multiples of $(a,b)$ also work for $s \neq 1$ and so we have that
    $\Lambda^\circ = \Lambda$. Since gcd$(a,b) = 1$, we have that $m = n$ if and only if
    $m = n = 0$ and so $A_p(\vphi)[m,n] = 0$ for every $(m,n) \in \Lambda^\circ$ except at
    $(0,0)$. Since $|\Lambda| = N$, we have that $(\vphi,\Lambda)$ is a tight frame with
    frame bound $N$. $\blacksquare$
    \end{prf}


The next example utilizes the square length Bj\"orck-Saffari sequence.
From Section \ref{slbsdpaf}, the DPAF of the square length Bj\"orck-Saffari sequence is
still very sparse, and the nonzero entries have regular structure as long as we take
the identity permutation $\sigma(m)= m$. We leverage this fact to prove \thref{sqlenframe} and
apply the same techniques as in \thref{chup4ap}.

\begin{prop}\thlabel{sqlenframe}
Let $c \in \C^N$ and let $\vphi \in \C^{N^2}$ be the square length Bj\"orck-Saffari sequence
generated by $c$ and $\sigma(h) = h$.
Let $K = \{0, a, \cdots, (bN'-1)a\}$ and $L = \{0, b, \cdots, (aN'-1)b\}$ with gcd$(a,b) = 1$ and
$N^2 = ab N'$. Then, the Gabor system $(\vphi, K \times L)$ is a tight Gabor frame with frame bound
$N^2 N'$.
\end{prop}

\begin{prf}
    Using \thref{72sol}, $\sigma(h) = h$, and the results of Section \ref{slbsdpaf},
    if $(rN'a,sN'b) \in N'K \times N'L$, then
    $A_p(\vphi)[rN'a,sN'b]$ can only be nonzero if
\begin{equation}\label{sqlenjan}
rN'a \equiv sN'b \mod N
\end{equation}
for some $r,s \in (\Z/N\Z)$. From the proof of \thref{chup4ap}, we know that
$N'K \cap N'L$ is generated by lcm$(aN',bN')$. We have that
lcm$(N'a,N'b) = N'$ lcm$(a,b) = abN' = N^2$, and we have that $N^2 \equiv 0 \mod N$.
Moreover, if $A_p(\vphi)[rN'a,sN'b] \neq 0$, then $r=s=0$. Appealing once again to
\thref{aptightframe}, we have that $(\vphi, K \times L)$ is a tight Gabor frame, and the
frame bound is $|\Lambda| = N^2 N'$. $\blacksquare$
\end{prf}


We now give an example utilizing \thref{aptightframe}, but with the opposite theme: the
discrete periodic amibguity function will be mostly nonzero except for the ``right'' spots.
That is, $A_p(\vphi)[m,n] \neq 0$ for nearly every $(m,n) \not\in \Lambda^\circ$ but
$A_p(\vphi)$ will still be $\Lambda^\circ$-sparse.

\begin{prop}
    Let $u \in \C^M$ be unimodular and satisfy
    \begin{equation}\label{uzac}
        \forall k \in (\Z/N\Z), k \neq 0, \text{ } 
        \langle \tau_{-k}u, u \rangle = N A_p(u)[k,0] = 0. 
    \end{equation}
    Furthermore, let $v\in \C^N$ be also unimodular and
    let $\vphi \in \C^{MN}$ be defined by
\[
\vphi = u \otimes v
\]
where $\otimes$ is the Kronecker product.
Assume gcd$(M,N) = 1$ and let  $\Lambda = \{0,M,\cdots,(N-1)M\} \times
\{0,N,\cdots,(M-1)N\}$. 
Then, $(\vphi, \Lambda)$ is a tight
frame with frame bound $MN$.
\end{prop}

\begin{prf}
We can write the $(rN + s)$-th term of $\vphi$ as
\[
    \vphi[rN+s] = u[r]v[s]
\]
where $r \in (\Z/M\Z)$ and $s \in (\Z/N\Z)$. We now compute the DPAF,
\[
    A_p(\vphi)[rN+s,\ell] = \frac{1}{MN} \sum_{j = 0}^{M-1} \sum_{k = 0}^{N-1}
    u\left[j+r+ \left\lfloor \frac{s+k}{N} \right\rfloor\right] v[s+k] \overline{u[j]v[k]}
    e^{-2 \pi i (jN+k) \ell /MN}
\]
\[
    = \frac{1}{N} \sum_{k = 0}^{N-1} v[s+k] \overline{v[k]} e^{-2 \pi i k \ell / MN}
    \frac{1}{M}
    \sum_{j=0}^{M-1} u\left[j + r+ \left\lfloor \frac{s+k}{N} \right\rfloor\right]
    \overline{u[j]} e^{-2 \pi i j \ell / M}
\]
\begin{equation}\label{krondpaf}
    = \frac{1}{N} \sum_{k=0}^{N-1} v[s+k] \overline{v[k]} e^{-2 \pi i k \ell / MN}
    A_p(u)\left[r+\left\lfloor\frac{s+k}{N}\right\rfloor, \ell\right].
\end{equation}
By \thref{72sol} we have that
$\Lambda^\circ = \Lambda$. Using (\ref{krondpaf}), we can see that for $(rN,\ell M) 
\in \Lambda^\circ \setminus \{(0,0)\}$ we have

\[
    A_p(\vphi)[rN,\ell M] = \frac{1}{N}\sum_{k=0}^{N-1} v[k]\overline{v[k]}
    e^{-2 \pi i k \ell / N} A_p(u)[r,\ell M]
    = \frac{A_p(u)[r,0]}{N}\sum_{k=0}^{N-1} |v[k]|^2
    e^{-2 \pi i k \ell / N} = 0
\]
since one of $r$ or $\ell$ is nonzero. Indeed, if $r$ is nonzero, then by (\ref{uzac}) we 
have a multiplier
of zero outside of the sum and if $r = 0$ but $\ell \neq 0$ then the sum will add up to zero since
$|v[k]|^2 = 1$ for every $k$. We now conclude by \thref{aptightframe} that $(\vphi, \Lambda)$ is
a tight frame with frame bound $MN$. $\blacksquare$
\end{prf}

\begin{figure}[h!]
\centering
\includegraphics[width=0.9\textwidth]{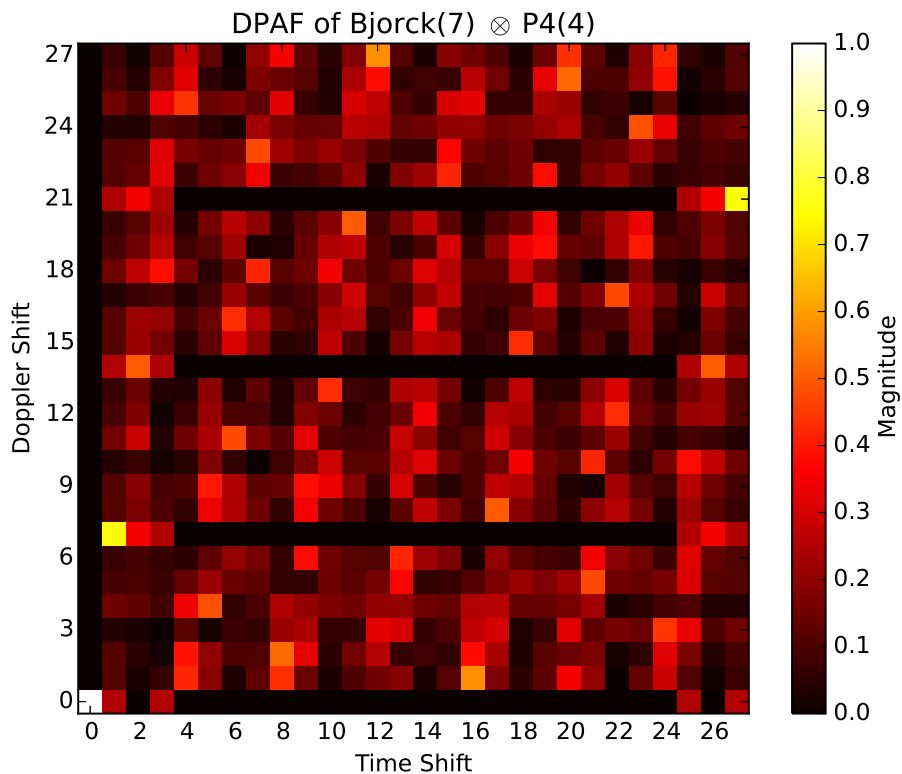}
\caption{DPAF of $\vphi = u \otimes v$ where $u$ is the length 7 Bj\"orck sequence
and $v$ is the length 4 P4 sequence.}
\label{kronaf}
\end{figure}

Although $A_p(\vphi)$ is $\Lambda^\circ$-sparse, 
(\ref{krondpaf}) implies that most of the entries for $A_p(\vphi)$ are nonzero. 
This is illustrated by Figure \ref{kronaf}. In Figure \ref{kronaf}, $u$ is the Bj\"orck sequence 
which is defined in Section \ref{cazacsec}. 
The definition of the Bj\"orck sequence and some of its relevant properties are defined in 
can also be found in \cite{BenKonRan2009}.
It should be noted that the Bj\"orck sequence is not the same as the square length
Bj\"orck-Saffari sequences defined in Section \ref{cazacdpaf}.

\begin{prop}
    Let $\vphi \in \C^{MN^2}$ be a Milewski sequence where the generating $v \in C^M$
    is either the Chu or P4 sequence, and let 
    $a,b,j,N'$ be such that gcd$(a,b) = 1$, $j \mid MN$, $N' = jN$, and $abN' = MN^2$. Let
    $\Lambda = K \times L$ where $K = \{0,a,\cdots,(bN'-1)a\}$, and $L = \{0, b, \cdots, (aN'-1)b 
    \}$. Then, $(\vphi, \Lambda)$ is a tight frame with frame bound $jMN^3$.
\end{prop}

\begin{prf}
By \thref{72sol}, we have that
$\Lambda^\circ = N'K \times N'L$. In particular, for every $(m,n) \in \Lambda^\circ$ we have
that $m \equiv n \equiv 0 \mod N$. Thus, using the second line of (\ref{milaf}), 
for $(m,n) \in \Lambda^\circ$
we can write $A_p(\vphi)[m,n]$ in the form
\[
    N A_p(\vphi)[m,n] = \sum_{j=0}^{N-1} e^{\theta[j]} A_p(v)[m',n']
\]
where $m' = m/N$ and $n' = n/N$. Since $v$ is the Chu or P4 sequence, we have that
$A_p(v)[m',n'] \neq 0$ if and only if $m' \equiv n' \mod M$. 
Furthermore, $\Lambda^\circ = \langle aN' \rangle \times \langle bN' \rangle
= \langle ajN \rangle \times \langle bjN \rangle$ and so we have that $m' \in \langle aj \rangle$ 
and $n' \in \langle bj \rangle$. Thus, we have that $m' \equiv n'$ if and only if $m',n' \in
\langle aj \rangle \cap \langle bj \rangle = \langle \text{lcm}(aj,bj) \rangle$.
Since gcd$(a,b) = 1$, we have that lcm$(aj,bj) = abj = MN$ since $abjN = abN' = MN^2$. 
In particular, we want to view $\langle aj \rangle$ and $\langle bj \rangle$ as being
subgroups of $(\Z/MN\Z)$ and in light of this we have that $m' \equiv n' \mod MN$ if and only 
if $m' \equiv n' \equiv 0 \mod MN$.
Thus,
we have that for $(m,n) \in \Lambda^\circ$, $A_p(\vphi) = 0$ unless $(m,n) = (0,0)$.
Using \thref{aptightframe}, we finally conclude that
$(\vphi, \Lambda)$ is a tight frame with frame bound $jMN^3$. $\blacksquare$
\end{prf}


\section{Gram matrices in terms of the discrete periodic ambiguity function}\label{gramdpaf}

The last three sections are devoted to an alternative method for showing when Gabor systems
are tight frames. The general framework is as follows: First, we explicitly compute the entries
of the Gram matrix by using the discrete periodic ambiguity function. Then, we show that the
first $N$ columns or rows happen to have disjoint supports. Last, we show that every other
column or row is a constant multiple of one of the first $N$ rows or columns and use this to
show that all $N$ nonzero eigenvalues are the same. This allows us to conclude that the Gabor
system is indeed a tight frame. We begin by defining the frame operator and Gram operator and
comparing the two.


\subsection{Gram Operator and Frame Operator} 
Let $\mathcal{H}$ be an $N$-dimensional Hilbert space and let $\mathcal{F} = \{v_i\}_{i=1}^M
\subseteq \mathcal{H}$ be
a frame for $\mathcal{H}$, where $M \geq N$. We define the \textit{analysis operator} 
$F: \mathcal{H} \to \R^M$ by
\begin{equation}\label{analysisop}
    \forall x \in \mathcal{H}, \text{ } F(x) = \{\langle x, v_i \rangle\}_{i=1}^M.
\end{equation}
The adjoint of the analysis operator, $F^*: \R^M \to \mathcal{H}$,
is called the \textit{synthesis operator} and is given by
\begin{equation}
    \forall \{c_i\}_{i=1}^M \in \R^M, \text{ } F^*(\{c_i\}_{i=1}^M) = \sum_{i=1}^M c_i v_i.
\end{equation}
Given the analysis and the synthesis operator, we can now define the frame operator
of the frame $\mathcal{F}$. The \textit{frame operator}, $S: \mathcal{H} \to \mathcal{H}$,
is defined by $S = F^* F$. We can write this explicitly as
\begin{equation}\label{frameop2}
    \forall x \in \mathcal{H}, \text{ } S(x) = \sum_{i=1}^M \langle x, v_i \rangle v_i.
\end{equation}
Note that $S$ is a self-adjoint operator. Indeed, $S^* = (F^* F)^* = F^* F^{**} = F^* F = S$.
We can define a second operator by reversing the order of the analysis and synthesis operator. 
That is, we apply the synthesis operator first
and the analysis operator second. This new operator $G: \R^M \to \R^M$ is called the
\textit{Gram operator} and is defined by $G = F F^*$. We can write this explicitly as
\begin{equation}\label{gramop}
    \forall \{c_i\}_{i=1}^M \in \R^M, \text{ } F(\{c_i\}_{i=1}^M) =
    \left\{ \left\langle \sum_{i=1}^M c_i v_i, v_j
    \right\rangle \right\}_{j=1}^M.
\end{equation}
(\ref{gramop}) is unwieldy and so it is usually more convenient to write the Gram operator in
matrix form. Once can see from (\ref{gramop}) that we can write the $(i,j)$-th entry of the
matrix form of $G$ as $G_{i,j} = \langle v_i, v_j \rangle$. A detailed exposition on finite
frames, the frame operator, and the Gram operator can be found in the first chapter of
\cite{Chr2008}, but we will use the following fact in Section \ref{secgramchu}.

\begin{thm}\thlabel{gramtightframes}
    Let $\mathcal{H}$ be an $N$-dimensional Hilbert space and let $\mathcal{F}
    = \{v_i\}_{i=1}^M$ be a frame for $\mathcal{H}$. $\mathcal{F}$ is a tight frame if and only
    if rank$(G) = N$ and every nonzero eigenvalue of $G$ is equal.
\end{thm}


\subsection{Gram Matrix of Gabor Systems}
In this section, we show that each entry of the Gram matrix of a Gabor system can be written in 
terms of the discrete periodic ambiguity function of the $\vphi \in \C^N \setminus \{0\}$ which
generates the system.
For all that follows let $\vphi \in \C^N \setminus \{0\}$ and let us write out the Gabor system as
$\mathcal{F} = \{e_{\ell_m}\tau_{k_m}\vphi : 
m \in {0, \cdots, M}\}$. 
Then, we can compute the $(m,n)$-th entry of the Gram matrix:
\[
G_{m,n} = \langle e_{\ell_m}\tau_{k_m}\vphi, \overline{e_{\ell_n}\tau_{k_n}\vphi} \rangle = 
\sum_{j = 0}^{N} e^{2 \pi i \ell_m j / N} 
\vphi[j-k_m] e^{-2 \pi i \ell_n j /N} \overline{\vphi[j-k_n]}
\]
\[
= \sum_{j=0}^N e^{-2 \pi i (\ell_n - \ell_m) j /N} \vphi[j-k_m] \overline{\vphi[j-k_n]} 
= \sum_{j=0}^N e^{-2 \pi i (\ell_n - \ell_m) (j-k_n+k_n) /N} \vphi[j-k_n+(k_n-k_m)] \overline{\vphi[j-k_n]} 
\]
\[
= e^{-2 \pi i k_n(\ell_n-\ell_m)/N} \sum_{j=0}^N \vphi[j+(k_n-k_m)]\overline{\vphi[j]}e^{-2 \pi i 
(\ell_n - \ell_m) j / N}
= N e^{-2 \pi i k_n(\ell_n-\ell_m)/N} A_p(\vphi)[k_n - k_m, \ell_n - \ell_m]
\]


\subsection{Gram Matrix for Chu Sequences}\label{gmtxchu}
In the Chu sequence case, we have that $G_{m,n} \neq 0$ if and only if $(\ell_n-\ell_m) 
\equiv (k_n-k_m) \mod N$. For convenience, let us
define $r_{mn}:\equiv k_n - k_m \equiv \ell_n - \ell_m \mod N$. Then, the Gram 
matrix for the Chu sequence is,
\begin{equation}
G_{mn} =
\begin{cases}
N e^{-2 \pi i [k_n r_{mn} - \frac{1}{2}(r_{mn}^2 - r_{mn})] / N}, & \text{if } (\ell_n- \ell_m) 
\equiv (k_n - k_m) \mod N \\
0, & \text{otherwise.} \label{chugram}
\end{cases}
\end{equation}


\subsection{Gram Matrices for P4 Sequences}\label{gmtxp4}
In the P4 sequence case, we also have that $G_{m,n} \neq 0$ if and only if $(\ell_n - 
\ell_m) \equiv (k_n - k_m) \mod N$. Again, 
let $r_{mn} :\equiv k_n-k_m \equiv \ell_n - \ell_m \mod N$. Then, the Gram 
matrix for the P4 sequence is,
\begin{equation}
G_{mn} = 
\begin{cases}
N (-1)^{r_{mn}} e^{-2 \pi i [k_n r_{mn} - \frac{1}{2}r_{mn}^2]/N}, & \text{if } (k_n - k_m) \equiv 
(\ell_n - \ell_m) \mod N \\
0, & \text{otherwise.} \label{p4gram}
\end{cases}
\end{equation}


\subsection{Gram Matrices for Wiener Sequences}\label{gmtxwiener}
In the Wiener sequence case, we have the following formula for the Gram matrices.
If $N$ is odd, then $G_{mn} \neq 0$ if and only if $2s(k_n - k_m) \equiv (\ell_n - \ell_m) \mod N$.
Let $r_{mn} :\equiv (k_n - k_m) \mod N$. Then we can write $G_{mn}$ as
\[
    G_{mn} = \begin{cases}
        N e^{-4 \pi i s k_n r_{mn}/N } e^{2 \pi i s r_{mn}^2 / N}, & \text{ if }
        2s(k_n - k_m) \equiv (\ell_n - \ell_m) \mod N \\
        0, & \text{ otherwise.}\\
    \end{cases}
\]
If $N$ is even, then then $G_{mn} \neq 0$ if and only if $s(k_n - k_m) \equiv (\ell_n - \ell_m)
\mod N$. Let $r_{mn}:\equiv (k_n - k_m) \mod N$. Then we can write $G_{mn}$ as
\[
    G_{mn} = \begin{cases}
        N e^{-2 \pi i s k_n r_{mn} / N} e^{\pi i s m^2 / N}, & \text{ if } 
        s(k_n - k_m) \equiv (\ell_n - \ell_m) \mod N \\
        0, & \text{ otherwise.} \\
    \end{cases}
\]
\section{Gram matrix method} \label{secgramchu}

In this section, we begin by applying the method outlined in Section \ref{gramdpaf} and apply it to 
the Chu and P4 sequence. We treat these two cases simultaneously since both have
the property that $A_p(\vphi)[m,n] = 0$ if $m \neq n$, and are nonzero when $m = n$. The
computation in \thref{chufnc} will be different for the P4 case, but the same ideas
can be applied to do the computation in the P4 case and achieve the same result. 
In the last part of this section we briefly discuss the Wiener sequence case.
The Wiener case has direct analogues of the results in the Chu and P4 cases and we will
only highlight the key differences in the proofs rather than reiterate all of the details.
Before proceeding, we will need a useful fact about the Gram matrix which can easily be derived 
from the singular value decomposition of the analysis operator. 

\begin{lemma}\thlabel{svdrank}
    Let $F$ be an $m \times n$ complex-valued matrix and let $G := FF^*$. Then, 
    rank($G$) = rank($F$).
\end{lemma}

\begin{prf}
First, we write $F$ in terms of its SVD: $F = UDV^*$ where $D$ is an $m \times n$ 
rectangular diagonal matrix and $U$ and $V$ are 
$m \times m$ and $n \times n$ unitary matrices, respectively. Note that $G = FF^* = UDV^*VD^*U^* 
= UDD^*U^* = UD^2U^*$. In particular, $D^2_{ij} =
|D_{ij}|^2$ and from this we get that $\text{rank}(G)=\text{rank}(F)$. $\blacksquare$
\end{prf}


For the following propostions, we shall use the following arrangement for the Gram matrix. 
Let $(\vphi, K \times L)$ be a Gabor system in $\C^N$. 
We shall iterate first by modulation, then iterate through translations.
In other words, the analysis operator, $F$, will be a $|K||L| \times N$ matrix where the $m$-th
row is given by $m = r|L| + s$, where $r \in \{0, \cdots, |K| - 1\}$,  $s = \{0, \cdots, |L| 
- 1 \}$, and the $m = r|L|+s$-th row corresponds to $e_{\ell_s} \tau_{k_r} \vphi$.


\begin{prop}\thlabel{abframe}
    Let $N = ab$ with gcd$(a,b) = 1$ and let $\vphi \in \C^N$ be either the 
    Chu or P4 sequence. Let $K = \{0,a,\cdots, (b-1)a\}$ and $L = \{0,b,\cdots,(a-1)b\}$.
    Then, the Gabor system $(\vphi, K \times L)$ is a tight frame with frame bound $N$.
\end{prop}

\begin{prf}
By construction, $|K| = b$, $|L| = a$ and so 
the Gabor system $(\vphi, K \times L)$ has $ab = N$ vectors. Since lcm$(a,b)=N$, we have that 
$K \cap L = \{0\}$.
From Section \ref{gmtxchu}, we have that
$G_{m,n} \neq 0$ if and only if $(\ell_n-\ell_m) = (k_n - k_m)$. However, by design of $K$ and $L$, 
we have that $(\ell_n - \ell_m) = jb$ and $(k_n - k_m)
= \tilde{j}a$ for some $j,\tilde{j}$. 
In particular, $(k_n - k_m) \in K$ and $(\ell_n - \ell_m) \in L$, and thus they are only equal
if they both belong to $K \cap L$.
From this we conclude that $G_{m,n} \neq 0$ if and only if 
$(\ell_n - \ell_m) = (k_n - k_m) = 0$, i.e. $j = \tilde{j} = 0$.
We conclude that the nonzero entries lie only on the diagonal of $G$.
Using formulas (\ref{chugram}) and (\ref{p4gram}), we see that for each $n \in \{0, \cdots, N-1\}$, 
$G_{n,n} = N$ and $G = N I\!d_{N}$.
Thus, the Gabor system $(\vphi, K \times L)$ is a tight frame with frame bound 
$N$. $\blacksquare$
\end{prf}


\begin{lemma}\thlabel{chufnc}
    Let $\vphi \in \C^N$ be either the Chu or P4 sequence and let $N = abN'$ where gcd$(a,b) = 1$.
    Suppose $G$ is the Gram matrix generated by the Gabor system $(\vphi, K \times L)$ where
    $K = \{0, a,\cdots,(bN'-1)a\}$ and $L = \{0,b,\cdots,(aN'-1)b\}$.
    Then, there exist functions $f,g: \N \cup \{0\} \to \C$ such that $G_{mn}/N = 
    f(n) g(m)$ wherever $G_{mn} \neq 0$.
\end{lemma}

\begin{prf}
We will only cover the case of the Chu sequence. 
The case of the P4 sequence follows by replacing $e^{-\pi i r_{mn} / N}$ with $(-1)^{r_{mn}}$ and 
carrying out the same computations.
If $G_{mn} \neq 0$, then we have
\begin{equation}
G_{mn}/N = e^{-2\pi i k_n r_{mn} / N} e^{-\pi i r_{mn}/N} e^{\pi i r_{mn}^2 / N} \label{gramchu}
\end{equation}
where 
\begin{equation}
r_{mn} = k_n - k_m = \ell_n - \ell_m = jab. \label{chucond} 
\end{equation}
Note that $k_n = c_nab + d_na$ and $k_m = c_mab + d_ma$ where $c_n,c_m \in \{0,\cdots,N'-1\}, 
d_n,d_m 
\in \{0,\cdots,b-1\}$. However, by (\ref{chucond}), $j = c_n - c_m$ and $d_n = d_m$.
Putting this back into (\ref{gramchu}) we have
\[
G_{mn}/N = e^{-2 \pi i (c_n ab + d_na)(c_n - c_m)ab/N}e^{-\pi i (c_n - c_m)ab/N} 
e^{\pi i (c_n - c_m)^2 a^2 b^2 / N}
\]
\[
= e^{-2 \pi i c_n ab (c_n - c_m)/N'} e^{-2 \pi i d_n a (c_n - c_m) / N'} e^{-\pi i (c_n - c_m)/N'} 
e^{\pi i ab (c_n^2 - 2c_nc_m + c_m^2)/N'}
\]
\[
= e^{-2 \pi i c_n^2 ab/N'}e^{2 \pi i c_n c_m ab / N'} e^{-2 \pi i d_n a (c_n - c_m) / N'} 
e^{-\pi i (c_n - c_m)/N'} e^{\pi i ab (c_n^2 + c_m^2)/N'}
e^{-2 \pi i ab c_n c_m / N'}
\]
\[
= e^{-\pi i c_n^2 ab/N'} e^{-2 \pi i d_n c_n a / N'} e^{-\pi i c_n/N'} e^{\pi i c_m^2 ab/N'} 
e^{2 \pi i d_m c_m a / N'} e^{\pi i c_m/N'}.
\]
Here we have used the crucial fact that $d_n = d_m$. Thus, we can write nonzero entries of 
$G_{mn}$ as $f(n)g(m)$ where
\[
    f(n) = e^{-\pi i c_n^2 ab/N'} e^{-2 \pi i d_n c_n a /N'}e^{-\pi i c_n / N'}
\]
and
\[
    g(m) = e^{\pi i c_m^2 ab/N'}e^{2 \pi i d_m c_m a / N'} e^{\pi i c_m / N'}. \blacksquare
\]
\end{prf}


\begin{lemma}\thlabel{chudisj}
    Let $\vphi \in \C^N$ be either the Chu or P4 sequence and let $N = abN'$ where gcd$(a,b) = 1$.
    Suppose $G$ is the Gram matrix generated by the Gabor system $(\vphi, K \times L)$ where
    $K = \{0, a,\cdots,(bN'-1)a\}$ and $L = \{0,b,\cdots,(aN'-1)b\}$.
    Then, the support of the rows (or columns) of $G$ either completely coincide or
    are completely disjoint.
\end{lemma}

\begin{prf}
Let us denote the $m$-th row of $G$ as $g_m$, and all
indices are implictly taken modulo $N$. We shall show that
if there is at least one index $n \in (\Z / N\Z)$ such that $g_m[n] = g_{m'}[n] \neq 0$, 
and $m \neq m'$, 
then supp$(g_m) = $ supp$(g_{m'})$.
Suppose that $g_m[n] = g_{m'}[n] \neq 0$.
Then, $k_n - k_m = \ell_n - \ell_m = jab$ and $k_n - k_{m'} = 
\ell_n - \ell_{m'} = j'ab$. We can rearrange these two equations to obtain the following
\begin{equation}
    k_m = k_n - jab \label{kmeqn}
\end{equation}
and
\begin{equation}
    k_n = k_{m'} + j'ab \label{kmpeqn}.
\end{equation}
Now suppose there is another $n' \in (\Z/N\Z)$ where $g_m[n'] \neq 0$. Then, we have
\begin{equation}
    k_{n'} - k_m = \tilde{j}ab  \label{kneqn}.
\end{equation}
Substiuting (\ref{kmeqn}) and (\ref{kmpeqn}) into (\ref{kneqn}), we have
\[
    \tilde{j}ab = k_{n'} - k_m = k_{n'} - k_n + jab = k_{n'} - k_{m'} - j'ab + jab,
\]
which can be rearranged to obtain
\[
    (j' + \tilde{j} - j)ab = k_{n'} - k_{m'}.
\]
Note that $(j'+\tilde{j}-j) \in (\Z/N\Z)$ and that these same computations can be done replacing $k$
with $\ell$. Thus, $g_{m'}[n'] \neq 0$ as well and the result is proved. $\blacksquare$
\end{prf}


\begin{rmk} 
Let gcd$(a,b)=1$, $N = abN'$, 
$K = \{0, a, \cdots, (bN'-1)a\}$, $L = \{0, b, \cdots, (aN'-1)b\}$, $\vphi \in \C^N$ be the
Chu and P4 sequence, and consider the system $(\vphi, K \times L)$. In light of
\thref{chufnc} and \thref{chudisj}, 
if two rows of the Gram matrix $G$ have supports which coincide,
they must be constant multiples of each other which has modulus 1.
Indeed, if $g_m$ and $g_{m'}$ have coinciding supports then for each $n$ where they are nonzero
we have
\[
    \frac{G_{mn}}{G_{m'n}} = \frac{Nf(n)g(m)}{Nf(n)g(m')} = \frac{g(m)}{g(m')}
\]
which also gives us a formula for finding the constant multiple, should we desire it.
This idea is illustrated with Figure \ref{picgram}, where two sets of rows with coinciding
supports are highlighted in red and blue, respectively.
\begin{figure}[h!]
\centering
\includegraphics[width=0.9\textwidth]{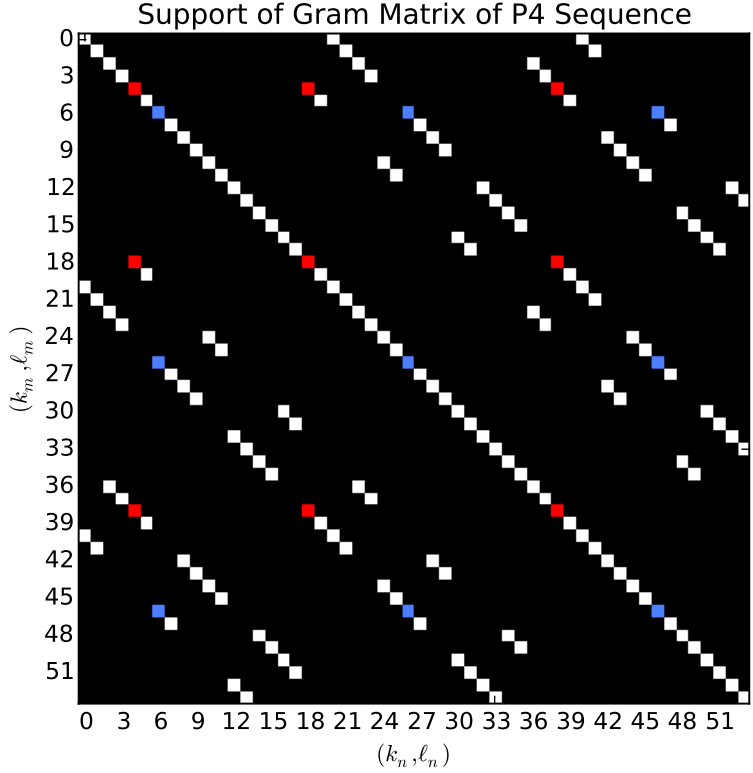}
\caption{Support of $G$ in the case $N$ = 18, $\vphi$ is the P4 sequence, $a = 2$, $b=3$,
        and $N' =3$.}
\label{picgram}
\end{figure}
\end{rmk}

\begin{thm}\thlabel{chuframe}
    Let $K = \{0,a,\cdots,(b-1)a\}$, $L = \{0,b,\cdots,(a-1)b\}$, and $N = abN'$ with
    gcd$(a,b) = 1$. Furthermore, let $\vphi \in \C^N$ be either the Chu or P4 sequence. Then,
    the Gabor system $(\vphi, K \times L)$ has the following properties:
    \begin{itemize}
        \item[(i)] rank$(G) = N$.
        \item[(ii)] The nonzero eigenvalues of $G$ are $NN'$.  
    \end{itemize}
    In particular, (i) and (ii) together imply that the Gabor system $(\vphi, K \times L)$ is a
    tight frame with frame bound $NN'$.
\end{thm}

\begin{prf}
(i) We shall show that the first $N$ columns
are disjoint, and then conclude that rank$(G) = N$. 
Note that $K$ and $L$ are subgroups of $(\Z/N\Z)$ and that $K \cap L = 
\{0, ab, 
2ab, \cdots, (N'-1)ab\}$. Moreover, 
$G_{mn} \neq 0$ if and only if $(\ell_n - \ell_m) \equiv (k_n - k_m) \mod N$. Since $K$ and
$L$ are subgroups, this can only happen if the subtractions lie in the intersections of the two
groups. That is, $G_{mn} \neq 0$ if and only if
\begin{equation}
(\ell_n - \ell_m) \equiv (k_n - k_m) \equiv jab \mod N \label{condi}
\end{equation}
where $j \in \{0, \cdots, N'-1\}$. Let
$m \in (\Z/N\Z)$. By $(\ref{condi})$, we have that for $G_{mn} \neq 0$ we must have that 
$k_n = (k_m + jab)$ for some $j$. Since $k_m 
\in K$, we have that 
\begin{equation}
    k_n = (j_m a + jab) = a (j_m + jb) \label{colform}
\end{equation}
for some $j_m \in \{0,\cdots,bN'-1\}$ and some $j \in \{0,\cdots,N'-1\}$.

By the ordering we used for the columns of $G$, we can write the index for column $n$ in terms of
$k_n$ and $\ell_n$ by
\begin{equation}
    n = (k_n/a) aN' + \ell_n/b = k_nN' + \ell_n/b.
\end{equation}
In particular, we would like $n \leq N$, and for that we need that $k_n < ab$. Looking at
(\ref{colform}), we need $(j_m + jb) < b$. There is exactly 1 such 
$j \in \{0, \cdots N'-1\}$ 
which can achieve this and it is obtained by setting $j = -\lfloor{j_m/b}\rfloor$. Thus, for each
row $m$, there is exactly 1 column $n \leq N$ where $G_{mn} \neq 0$, and therefore the first $N$
columns of $G$ are linearly independent. Thus, we conclude $\text{rank}(G)\geq N$. By
\thref{svdrank}, we know that rank($G$) = rank($L$). Since $L$ is an $M \times N$ matrix, we have
that $\text{rank}(G)\leq N$ and we have that $\text{rank}(G) = N$. We finally conclude
the Gabor system in question forms a frame.

(ii) Let $g_n$ be the $n$-th column of $G$, with $n \leq N$.
We wish to show that $Gg_n = NN' g_n$.
Note that $Gg_n[m]$ is given by the inner product of the $m$-th row and the $n$-th column of $G$.
Furthermore, since $G$ is self-adjoint, the $n$-th column is also the conjugate of the $n$-th row.
If $g_n[m] \neq 0$, then $G_{mn} \neq 0$ and $G_{nn} \neq 0$. Thus, by \thref{chudisj} rows $m$ 
and $n$ have supports that coincide.
By \thref{chufnc}, $G_{m(\cdot)} = C_m g_n^*$, where $|C_m| = 1$.
Thus, $Gg_n[m] = C_m \|g_n\|_2^2 = N^2N' C_m$.
It is easily computed that $G_{nn} = N$, so by \thref{chufnc}, we have that $g_n[m] = NC_m$. In 
particular, we get
$Gg_n[m] = (NN')(NC_m)$.
Thus, the first $N$ columns of $G$ are eigenvectors of $G$ and they all have eigenvalue $NN'$.
It now follows that the system is a tight frame with frame bound $NN'$. $\blacksquare$
\end{prf}


\begin{rmk} 
In general, if $\gcd(a,b) > 1$, then the above result will not hold. Let 
$\vphi \in \C^4$ be the P4 sequence. That is,
$\vphi = (1, -\sqrt{2}/2 - i \sqrt{2}/2, -1, -\sqrt{2}/2 - i \sqrt{2}/2)$. Let $K = \{0,2\}$ and 
$L = \{0,2\}$. Note that this would give
$a = b = \gcd(a,b) = 2$. The Gabor system $(\vphi, K\times L)$ is given by
\[
(\vphi,K\times L) = \{\pi(0,0)\vphi, \pi(0,2)\vphi, \pi(2,0)\vphi, \pi(2,2)\vphi\} 
\]
\[
= \left\{
\begin{pmatrix}1\\-\frac{\sqrt{2}}{2}-\frac{\sqrt{2}}{2}i\\-1\\
-\frac{\sqrt{2}}{2}-\frac{\sqrt{2}}{2}i\end{pmatrix},
\begin{pmatrix}1\\\frac{\sqrt{2}}{2}+\frac{\sqrt{2}}{2}i\\-1\\
\frac{\sqrt{2}}{2}+\frac{\sqrt{2}}{2}i\end{pmatrix},
\begin{pmatrix}-1\\-\frac{\sqrt{2}}{2}-\frac{\sqrt{2}}{2}i\\1\\
-\frac{\sqrt{2}}{2}-\frac{\sqrt{2}}{2}i\end{pmatrix},
\begin{pmatrix}-1\\\frac{\sqrt{2}}{2}+\frac{\sqrt{2}}{2}i\\1\\
\frac{\sqrt{2}}{2}+\frac{\sqrt{2}}{2}i\end{pmatrix},
\right\}.
\]
Note that the first and fourth vectors are multiples of each other, as well as the second and
third vectors. Specifically, 
\[
\pi(0,0)\vphi = -\pi(2,2)\vphi \hs \text{ and } \hs \pi(0,2)\vphi = -\pi(2,0)\vphi.
\]
We conclude from this that the dimension of the span of the Gabor system $(\vphi, K \times L)$ is 
only 2, and the Gabor system in question is not a frame.
\end{rmk}


We close this section with brief mentions about the proofs in the Wiener sequence case.
To simplify further, we only mention the odd length case, but one can easily make the highlighted
adjustments in the even case as well. The emphasis here is that the results and proofs in the
Wiener cases are essentially the same.


\begin{lemma}\thlabel{wienerfnc}
    Let $\vphi \in \C^N$ be a Wiener sequence of odd length and let $N = abN'$ where gcd$(a,b) = 
    1$. Suppose $G$ is the Gram matrix generated by the Gabor system $(\vphi, K \times L)$,
    where $K = \{0, a, \cdots, (bN'-1)a\}$ and $L = \{0, b, \cdots, (aN'-1) b\}$. 
    Then, there exist functions $f,g: \N\cup\{0\} \to \C$ such that $G_{mn}/N = f(n)g(m)$
    wherever $G_{mn} \neq 0$.
\end{lemma}

    To prove \thref{wienerfnc},
    the same technique used in \thref{chufnc} of writing $k_n = c_n ab + d_n a$ and 
    $k_m = c_m ab + d_m a$ works here as well and the result follows from the same type of
    computations used in \thref{chufnc}.


\begin{lemma}\thlabel{wienerdisj}
    Let $\vphi \in \C^N$ be a Wiener sequence of odd length and let $N = abN'$ where gcd$(a,b)=1$.
    Suppose $G$ is the Gram matrix generated by $(\vphi, K \times L)$, where
    $K = \{0,a,\cdots,(bN'-1)a\}$ and $L = \{0,b,\cdots,(aN'-a)b\}$.
    Then, the support of the rows (and columns) of $G$ either completely coincide or are completely 
    disjoint.
\end{lemma}

To prove \thref{wienerdisj}, one needs to replace $jab$, $j'ab$, $\tilde{j}ab$ in
\thref{chudisj} with $2sjab$,
$2sj'ab$, and $2s\tilde{j}ab$, and the same result will hold.


\begin{thm}\thlabel{wienerframe}
    Let $K = \{0,a,\cdots,(b-1)a\}$, $L = \{0,b,\cdots,(a-1)b\}$, and $N = abN'$ with
    gcd$(a,b) = 1$. Furthermore, let $\vphi \in \C^N$ be a Wiener sequence of odd length. Then,
    the Gabor system $(\vphi, K \times L)$ has the following properties:
    \begin{itemize}
        \item[(i)] rank$(G) = N$.
        \item[(ii)] The nonzero eigenvalues of $G$ are $NN'$.
    \end{itemize}
    In particular, (i) and (ii) combined imply that the Gabor system $(\vphi, K \times L)$ is a tight
    frame with bound $NN'$.
\end{thm}

As with the modification used in \thref{wienerdisj}, one only needs to change any instance of
$jab$ in \thref{chuframe} with $2sjab$ and apply the appropriate computations to prove
\thref{wienerframe}.

\bibliographystyle{plain}
\bibliography{cazacframe}

\end{document}